\documentclass[11pt]{article}

\usepackage[pdftex]{graphicx}

\usepackage{textcomp}
\usepackage{amsmath}
\usepackage{amsfonts}
\usepackage{epsfig}

\renewcommand{\baselinestretch} {1.2}

\makeatletter 
\def\singlespace{\def\baselinestretch{1}\@normalsize}

\title{{\sc Bootstrap-Based $K$-Sample Testing   For Functional Data}}

\author{ {Efstathios ~{\sc PAPARODITIS}\footnote{
          Email: \texttt{stathisp@ucy.ac.cy}}} \; and \;  {Theofanis ~{\sc SAPATINAS}\footnote{Email: \texttt{fanis@ucy.ac.cy}}} \\
          Department of Mathematics and Statistics, University of Cyprus, \\
          P.O. Box 20537, CY 1678 Nicosia, CYPRUS.
}

\newcommand{\reals}{\ensuremath{{\mathbb R}}}
\newcommand{\RR}{\reals}

\newcommand{\EE}{\ensuremath{{\mathbb E}}}
\newcommand{\PP}{\ensuremath{{\mathbb P}}}

\newcommand{\la}{\mbox{$\langle$}}
\newcommand{\ra}{\mbox{$\rangle$}}
\newtheorem{theorem}{Theorem}[section]
\newtheorem{lemma}{Lemma}[section]
\newtheorem{remark}{Remark}[section]

\numberwithin{figure}{section}
\numberwithin{table}{section}

\begin{document}

\DeclareGraphicsExtensions{.pdf, .jpg}

\maketitle

\begin{abstract}
We investigate properties of a bootstrap-based methodology for testing hypotheses about equality of  certain characteristics of the distributions between different populations in the context of  functional data. The suggested testing methodology is simple and easy to implement. It  bootstraps the original functional dataset in such a way that the null hypothesis of interest is satisfied and it can be potentially applied to a wide range of testing problems and test statistics of interest. Furthermore,  it can be utilized to the case where more than two populations of functional data are considered.  We illustrate the bootstrap procedure by considering the important problems of testing the equality of  mean functions or the equality of covariance functions (resp. covariance operators) between two populations.  Theoretical results that justify the validity of the suggested bootstrap-based procedure are established. Furthermore, simulation results demonstrate very good size and power performances in finite sample situations, including the case of testing problems and/or  sample sizes  where  asymptotic considerations  do not lead  to satisfactory approximations. A real-life dataset analyzed in the literature is also examined. 

\medskip
\noindent
{\em Some key words:} {\sc 
Bootstrap; 
Covariance Function, 
Functional Data;  Functional Principal Components; 
Karhunen-L\`oeve Expansion; Mean Function; $K$-sample problem.}
\end{abstract}

\section{{\sc INTRODUCTION}}
\label{sec:intro}

Functional data are routinely collected in many fields of research; see, e.g., Bosq (2000), Ramsay \& Silverman (2002, 2005), Ferraty \& Vieu (2006), Ramsay {\em et al.} (2009) and Horv\'ath \& Kokoszka (2012).  They are usually recorded at the same, often equally spaced, time points, and with the same high sampling rate per subject of interest, a common feature of modern recording equipments.  The estimation of individual curves (functions) from noisy data and the characterization of  homogeneity and of patterns of variability among curves are main concerns of functional data analysis; see, e.g.,  Rice (2004). When working with more than one population (group), the equality of certain characteristics of the distributions between the populations, like their mean functions or their covariance functions (resp. covariance operators) is  an interesting and widely discussed problem in the literature; see, e.g., Benko {\em et al.} (2009), Panaretos {\em et al.} (2010), Zhang {\em et al.} (2010), Fremdt {\em et al.} (2012), Horv\'ath \& Kokoszka (2012), Kraus \& Panaretos (2012),  Horv\'ath {\em et al.} (2013), Fremdt {\em et al.} (2013) and Boente {\em et al.} (2014). 

\medskip

For instance, Benko {\em et al.} (2009) and  
 Horv\'ath \& Kokoszka (2012,  Chapter 5) have developed asymptotic functional testing procedures for the equality of two mean functions. For the more involved problem of testing the equality of covariance functions, Panaretos {\em et al.} (2010) and Fremdt {\em et al.} (2012) have developed corresponding  testing procedures in the two-sample problem. Critical points of these testing procedures are typically obtained using  asymptotic approximations of the distributions of the  test statistics used under validity of the null hypothesis. In this context, the main tools utilized are the functional principal components (FPC's) and the associated Karhunen-Lo\`eve expansion (KLE); see, e.g., Reiss \& Ogden (2007), Gervini (2008), Yao \& M\"uller (2010), Gabrys {\em et al.} (2010) and Fremdt {\em et al.} (2013). 


\medskip

For testing the equality of  two covariance functions, Panaretos {\em et al.} (2010) have derived a functional testing procedure  under the assumption of Gaussianity, while Fremdt {\em et al.} (2013)  have extended such a functional testing procedure to the non-Gaussian case.  Clearly,  and due to the complicated statistical functionals involved, the efficacy of these functional testing procedures heavily rely on the accuracy of the obtained asymptotic approximations of the distributions of the test statistics considered  under the null hypothesis.  Simulation studies, however,  suggest that the quality of some asymptotic approximations is questionable. This is not only true for  small or moderate sample sizes, that are of paramount importance in practical applications, but also in situations where the assumptions under which the asymptotic results  have been derived (e.g., Gaussianity) are not satisfied in practice; see, e.g., Fremdt  {\em et al.} (2013, Section 4).   

\medskip

To improve such asymptotic approximations, bootstrap-based testing  inference for  functional data have been considered by some authors in the literature.  For instance, Benko {\em et al.} (2009) have considered, among other things, testing the equality of mean functions in the two-sample problem and have applied a bootstrap procedure to obtain  critical values of the test statistics used. This bootstrap procedure resamples the original set of functional observations themselves without imposing the null hypothesis and, therefore, its  validity  rely  on the particular test statistic used. That is, the  bootstrap approach used does not generate  functional pseudo-observations that satisfy the null hypothesis of interest.  Therefore, it is not clear if  this  procedure  can be applied  to other  test statistics,  to  different testing problems  or  to  the case where  more than two populations of functional observations are compared.   Similarly, and for the case of comparing the mean functions of two populations of functional observations, Zhang {\em et al.} (2010) have considered  a bootstrap  procedure that generates functional pseudo-observations which do not satisfy the null hypothesis. Thus, the validity of this approach depends on the specific test statistic used.   A different idea for improving  asymptotic approximations has been used by  Boente {\em et al.} (2014) in the context of testing the equality of  several covariance functions, by applying a  bootstrap procedure in order to calibrate the critical values of the test statistic  used. Again, this bootstrap approach is taylor made for the particular test statistic considered and does involve any resampling of the  functional observations themselves.  Finally, permutation  tests for equality of covariance operators applied to different distance measures between two covariance functions  have been  considered by Pigoli {\em et al.} (2014).

\medskip

We investigate properties of an alternative and general bootstrap-based testing methodology for functional data, which is potentially applicable for different testing problems, different test statistics and for more than two populations.  Among other things, the bootstrap-based procedure proposed can be applied  to the important problem of comparing the mean functions or the covariance functions  between several populations. The basic idea behind this testing methodology is to bootstrap the observed functional data set in such a way that the obtained functional pseudo-observations  satisfy the   null hypothesis of interest.  This requirement leads to a particular bootstrap scheme that automatically generates pseudo-functional observations with  identical  mean functions (when testing the equality of mean functions) or identical covariance functions (when testing the equality of covariance functions) among the different populations. This common mean function is the estimated pooled mean function (when testing the equality of mean functions) and the common covariance function is the estimated pooled covariance function (when testing the equality of covariance functions) of the observed functional data. A given test statistic of interest is then calculated using the  bootstrap functional pseudo-observations and its distribution is evaluated by means of Monte Carlo simulations. As an example, we show that this bootstrap-based functional testing procedure consistently estimates the distribution of the test statistics under the null hypothesis, proposed by Benko {\em et al.} (2009) and Horv\'ath \& Kokoszka (2012,  Chapter 5), for the problem of testing the equality of two mean functions, and by Panaretos {\em et al.} (2011), Fremdt {\em et al.} (2013) and Boente {\em et al.} (2014), for the problem of testing the equality of two covariance functions. The theoretically established asymptotic validity of  the suggested bootstrap-based functional testing methodology is further gauged by extensive simulations and coincides with accurate approximations of the distributions of interest in finite sample situations. These accurate approximations lead to a very good size and power behavior of the test statistics considered.

%

\medskip

The paper is organized as follows. In Section \ref{sec:bft}, we first present   the suggested  bootstrap-based  methodology applied to  the problem of testing equality 
of  the covariance functions in the functional set-up. We then extend the discussion to the problem of testing equality of  the mean functions between several groups.  Some pertinent observations through remarks on the testing methodology considered are also included. In Section \ref{sec:theor}, we provide  theoretical results which justify  the validity of the suggested bootstrap-based  testing methodology applied  to  some test statistics recently considered   in the literature.  In Section \ref{sec:num}, we evaluate the finite sample behavior of the proposed  bootstrap-based testing procedures by means of several  simulations and  compare the results obtained with those  based on  classical asymptotic approximations.    An application to a real-life dataset is also presented. Some concluding remarks are made in Section \ref{sec:conc}. Finally, auxiliary results and proofs of the main results are compiled in the Appendix.

\section{\sc Bootstrap-based Functional Testing Methodology}
\label{sec:bft}
\subsection{\sc Model and Assumptions}
\label{subsubsec:m-a-m}
We work with functional data in the form of random functions $X:=X(t):=X(\omega, t)$, defined on a probability space $(\Omega, {\cal A}, \PP)$ with values in the separable Hilbert-space ${\cal H}=L^2({\cal I})$, the space of squared-integrable $\RR$-valued functions on the compact interval ${\cal I}=[0,1]$. We denote by $\mu(t):=\EE[X(t)]$, (for almost all) $t \in {\cal I}$, the mean function of $X$, i.e., the unique function $\mu \in L^2({\cal I})$ such that $\EE<X,x>=<\mu, x>$, $x \in L^2({\cal I})$.  We also denote by $C(t,s):=Cov[X(t),X(s)]:=\EE[(X(t)-\mu(t))(X(s)-\mu(s))]$, $t,s \in {\cal I}$,  the covariance function (kernel) of $X$, and by ${\cal C}(f)=\EE[\la X-\mu,f \ra (X-\mu)]$, for $f \in L^2({\cal I})$, the covariance operator of $X$. It is easily seen that ${\cal C}(f)(t)=\int_{{\cal I}}C(t,s)f(s)ds$, i.e., ${\cal C}$ is an integral operator with kernel $C$; note that $C$ is a Hilbert-Schmidt operator provided that $\EE||X||^2<\infty$. Throughout the paper we assume that
$\la f, g \ra = \int_{{\cal I}}f(t)g(t)dt$,  $\|f\|^2=\la f, f \ra$,
and that all functions considered are elements of the separable Hilbert-space $L^2({\cal I})$.
Finally, the operator $u \otimes v: L^2 \mapsto L^2$ is defined as $(u \otimes v)w=\langle v,w\rangle u$, $u,v \in L^2$, and we denote by $ \|{\mathcal C}\|_S$ the
 Hilbert-Schmidt norm of the covariance operator $ {\mathcal C}$.

\medskip

Throughout the paper, it is also assumed that we have available a collection of random functions satisfying 
\begin{equation}
\label{eq:model}
X_{i,j}(t)=\mu_i(t)+\epsilon_{i,j}(t), \quad i=1,2,\ldots,K, \;\; j=1,2,\ldots,n_i, \;\;  t \in {\cal I},
\end{equation}
where $K \;(2 \leq K < \infty)$ denotes the number of populations (groups), $n_i$ denotes the number of observations in the $i$-th population and $N=\sum_{i=1}^{K}n_i$ denotes the total number of observations. We also assume that the $K$ populations are independent and, for each $i \in \{1,2,\ldots,K\}$ and $j=1,2,\ldots,n_i$, the $\epsilon_{i,j}$ are independent and identical distributed random elements with $\EE[\epsilon_{i,j}(t)]=0$, $t \in {\cal I}$, and  $\EE\|\epsilon_{i,j} \|^4 <\infty$. 

\medskip

Denote by $(\lambda_k, \varphi_k)$, $k=1,2,\ldots,$ the eigenvalues/eigenfunctions of the covariance operator ${\cal C}$, i.e.
\[  
\lambda_k \varphi_k(t) = {\cal C}(\varphi_k)(t) :=\int_{\cal I} C(t,s) \varphi_k(s)ds, \ \ \ t\in {\cal I}, \ \ \ k=1,2,\ldots \ .
\]
Throughout the paper it is assumed that $\lambda_1 > \lambda_2 >  \cdots >\lambda_p > \lambda_{p+1}$,
i.e., there exists at least $p$ distinct (positive) eigenvalues of the covariance operator ${\cal C}$.

\subsection{\sc Testing the Equality of Covariance Functions}
\label{subset:cov}
In this section, we describe the suggested bootstrap-based functional testing methodology for testing    the equality of  covariance functions (resp. covariance operators) for a (finite) number of populations. Since  testing the equality of  covariance functions is equivalent to testing the equality of  covariance operators, as in Fremdt {\em et al.} (2012), we confine our attention to the former test.

\medskip

Let ${\bf X}_N=\{X_{i,j}(t), \, i \in \{1,2,\ldots,K\},  \, j=1,2,\ldots,n_i, \, t \in {\cal I}\}$ be  the observed  collection of random functions satisfying (\ref{eq:model}). Let $C_i(t,s)$, $t,s \in {\cal I}$, be the covariance functions in the $i$-th population, i.e., for each $i \in \{1,2,\ldots,K\}$, $C_i(t,s):=Cov[X_{i,j}(t),X_{i,j}(s)]:=\EE[(X_{i,j}(t)-\mu_i(t))(X_{i,j}(s)-\mu_i(s))]$,  where $\mu_i(t):=\EE[X_{i,j}(t)]$, $j=1,2,\ldots,n_i$. Our aim is  to test the null hypothesis 
\begin{equation}
\label{null-cov}
H_0:\; C_1=C_2=\ldots =C_K
\end{equation}
versus the alternative hypothesis
\begin{equation}
\label{alt-cov}
H_1:\; \exists \; (k,l)\in\{1,2,\ldots,K\} \;\; \text{with}\;\; k \neq l\;\; \text{such that}\;\; C_k \neq C_l .
\end{equation}

Notice that  the equality in the null hypothesis (\ref{null-cov}) is in the space $(L^2({\cal I}\times {\cal I}), \, \| \cdot \|)$, i.e., $C_k=C_l$, for any pair of indices $(k,l)\in\{1,2,\ldots,K\}$, with $k \neq l$, means that  $ \|C_k-C_l\| =0$, 
and the alternative hypothesis (\ref{alt-cov}) means that 
$ \|C_k-C_l\| >0$.

\subsubsection{\sc The Bootstrap-based Testing Procedure}
\label{subsec.bootalg}

 Let $T_N$ be a  given test statistic of interest  for testing hypothesis (\ref{null-cov})
which is  based on the functional observations $ {\bf X}_N$.
 Assume, without loss of generality, that $ T_N$  rejects the null hypothesis $H_0$ when $T_N >d_{N,\alpha}$,  where  for $ \alpha \in (0,1)$, 
 $d_{N,\alpha}$ denotes  the  critical value of this test. The bootstrap-based functional testing procedure for testing  hypotheses (\ref{null-cov})-(\ref{alt-cov}) can be described as follows:

\vspace*{0.3cm}

\renewcommand{\baselinestretch}{0.7}
\small\normalsize 

\medskip

\noindent {\bf Step 1:} \  First calculate the sample mean functions in each population 
$$\overline{X}_{i,n_{i}}(t)=\frac{1}{n_i}\sum_{j=1}^{n_i} X_{i,j}(t), \quad t \in {\cal I},\quad i \in \{1,2,\ldots,K\}.$$
\medskip
\noindent {\bf Step 2:} \ Calculate the residual functions in each population, i.e., for each $i \in \{1,2,\ldots,K\}$,
$$
\hat{\epsilon}_{i,j}(t)=X_{i,j}(t)-\overline{X}_{i,n_{i}}(t), \quad t \in {\cal I}, \quad j=1,2,\ldots,n_i.
$$ 
\medskip
\noindent {\bf Step 3:} \ Generate bootstrap functional pseudo-observations  $X^*_{i,j}(t)$, $ t \in {\cal I}$, $i \in \{1,2,\ldots,K\}$, $j=1,2,\ldots,n_i$, according to
\begin{equation} \label{eq.xboot}
X^*_{i,j}(t)=\overline{X}_{i,n_i}(t) + \epsilon^*_{i,j}(t), \quad  t \in {\cal I},
\end{equation}
where 
$$
\epsilon^*_{i,j}(t)=\hat{\epsilon}_{I,J}(t), \; t \in {\cal I}, $$
and  $ (I,J) $ is the following   pair of   random variables.  
The random variable $ I$  takes values in the set $\{1,2,\ldots, K\}$ with probability  $ P(I=i)=n_i/N$ for $ i \in \{1,2, \ldots, K\}$, and, given  $I=i$,  the random variable $ J$ has the discrete uniform distribution 
in the set  $ \{1,2, \ldots, n_i\}$, i.e., $ P(J=j \mid I=i)= 1/n_i$ for $i \in \{1,2, \ldots, K\}$, $ j=1,2, \ldots, n_i$.
\medskip

\noindent {\bf Step 4:} \ Let $ T^\ast_N$ be the same statistic as  $T_N$ but  calculated using  the bootstrap functional pseudo-observations  $X^*_{i,j}$, $i \in \{1,2,\ldots,K\}$; $j=1,2,\ldots,n_i$. Denote by $ D^\ast_{N,T}$ the distribution function of $ T^\ast_N$ given the functional observations $ {\bf X}_N$.
\medskip

\noindent {\bf Step 5:} \ For any given $\alpha \in (0,1)$, reject the null hypothesis  $H_0$ if and only if 
\[ T_N \  > \  d^{\ast}_{N,\alpha},\]
 where $d^\ast_{N,\alpha}$ denotes   the $\alpha$-quantile of $D^\ast_{N,T}$, i.e., $D^\ast_{N,T}(d^\ast_{N,\alpha})=1-\alpha$.

\medskip
\renewcommand{\baselinestretch}{1.3}
\small\normalsize
Notice that   the distribution  $D^\ast_{N,T} $ can be   evaluated by means of Monte-Carlo.


\vspace*{0.3cm}

Clearly, and  since  the random  functions $\epsilon^\ast_{i,j}(t)$ are generated  independently  from each other, 
for any two different  pairs  of indices, say $(i_1,j_1)$ and $ (i_2,j_2)$,  the corresponding 
bootstrap functional pseudo-observations  $ X^\ast_{i_1,j_1}(t)$ and $ X^\ast_{i_2,j_2}(t)$ are independent. 
Furthermore,  observe  that the random selection of  the error function $\epsilon^\ast_{i,j}(t)$ in Step 4  of the above bootstrap algorithm, is equivalent to selecting  $\epsilon^\ast_{i,j}(t)$
randomly  with  probability $1/N$  from the entire set of   available and estimated residual  functions  $ \{ \widehat{\epsilon}_{r,s}(t):  r=1,2, \ldots, K \ \mbox{and} \ s=1,2, \ldots, n_r\}$. Hence,  conditional on the observed functional data ${\bf X}_N$, the functional pseudo-observations 
$ X_{i,j}^\ast(t)$ have the following first and second order properties:
$$
 \EE[X_{i,j}^\ast(t)]= \overline{X}_{i,n_i}(t) + \EE[\epsilon^*_{i,j}(t)] 
 = \overline{X}_{i,n_i}(t) + \frac{1}{N}\sum_{i=1}^{K}\sum_{j=1}^{n_i} \widehat{\epsilon}_{i,j}(t)= \overline{X}_{i,n_i}(t), \quad t \in {\cal I},
$$  
since     $ \sum_{j=1}^{n_i} \widehat{\epsilon}_{i,j}(t)= \sum_{j=1}^{n_i} (X_{i,j}(t)-\overline{X}_{i,n_i}(t))=0$, within each population $i \in \{1,2, \ldots, K\}$. \\Moreover, 
\begin{align*} 
Cov[X_{i,j}^\ast(t), X_{i,j}^\ast(s)]  & =  
 \EE[\epsilon_{i,j}^\ast(t)\epsilon_{i,j}^\ast(s)] = \frac{1}{N} \sum_{i=1}^{K}\sum_{j=1}^{n_i}\widehat{\epsilon}_{i,j}(t)\widehat{\epsilon}_{i,j}(s)\\
& = \sum_{i=1}^{K}\frac{n_i}{N} \frac{1}{n_i}\sum_{j=1}^{n_i}(X_{i,j}(t) -  \overline{X}_{i,n_i}(t))(X_{i,j}(s) - \overline{X}_{i,n_i}(s) ) \\
& = \sum_{i=1}^{K} \frac{n_i}{N} \widehat{C}_{i,n_i}(t,s) = \widehat{C}_N(t,s), \quad t,s \in {\cal I},
\end{align*}
where 
$$ \widehat{C}_{i,n_i}(t,s)=\frac{1}{n_i}\sum_{j=1}^{n_i}(X_{i,j}(t) -  \overline{X}_{i,n_i}(t))(X_{i,j}(s) - \overline{X}_{i,n_i}(s)), \quad t, s \in {\cal I}, 
$$ 
is the sample estimator of the covariance function $ C_{i}(t,s)$ for the $i$-th population and  $ \widehat{C}_N(t,s)$ is the corresponding pooled covariance function estimator.\\

Thus,  and conditional on the observed functional data $ {\bf X}_N$, the bootstrap generated functional pseudo-observations $ X^\ast_{i,j}(t)$ have, within each population 
$ i \in \{1,2, \ldots, K\}$, the same mean function $ \overline{X}_{i,n_{i}}(t)$, which may be different for different populations. Furthermore,  the covariance function in each population is identical  and equal to the pooled sample covariance function $ \widehat{C}_{N}(t,s)$.  That is, the functional pseudo-observations $ X_{i,j}^\ast(t)$, satisfy the  null hypothesis (\ref{null-cov}). This basic property of  the $ X^\ast_{i,j}(t)$'s allows us to use these bootstrap observations to evaluate the 
distribution  of some test statistic $ T_N$ of interest  under the   null hypothesis. This is achieved by 
using  the distribution of $ T_N^\ast$ as an  estimator of  the distribution of $ T_N$,  where $ T^\ast_N$ is the same statistical functional as $ T_N$  calculated using the  bootstrap functional   pseudo-observations $ {\bf X}^\ast_N=\{X_{i,j}^\ast(t),\, i=1,2, \ldots, K,\, j=1,2, \ldots, n_K, \ t \in {\cal I} \}$.   
Since, as we have seen,  the set of pseudo-observations used to calculate $T_N^\ast$ satisfy the null hypothesis,  we expect that the distribution of the pseudo random variable $ T_N^\ast$   will   mimic  correctly the distribution of $T_N$ under the null. In the next section we show that this is indeed true for two particular test statistics proposed in the literature. However, since our bootstrap methodology is not designed or tailor made for any particular test statistic,   its  range of validity  is not restricted to these 
two particular test statistics.

\subsection{\sc Testing the Equality of Mean Functions}
  \label{re.adapt}

We assume again that we have available a collection of curves ${\bf X}_N$,
satisfying  (\ref{eq:model}). Recall that $\mu_i(t)$, $t\in {\cal I}$ denote the mean functions of the curves  in the $i$-th population, i.e., for each $i \in \{1,2,\ldots,K\}$, $\mu_i(t):=\EE[(X_{i,j}(t)]$, $j=1,2,\ldots,n_i$. 

\medskip
The basic idea used in the bootstrap resampling algorithm of Section \ref{subsec.bootalg} enables  its adaption/modification   to deal with different testing problems  related to the comparison of $K$ populations of functional observations. For instance, suppose that we are interested in testing the null hypothesis 
that the $K$ populations have identical mean functions, i.e., 
\begin{equation} \label{null-mean}
H_0: \  \mu_1=\mu_2= \cdots = \mu_K
\end{equation}
versus the alternative hypothesis
\begin{equation}
\label{alt-mean}
H_1:\; \exists \; (k,l)\in\{1,2,\ldots,K\} \;\; \text{with}\;\; k \neq l\;\; \text{such that}\;\; \mu_k \neq \mu_l.
\end{equation}

As in the previous section,   equality in the null hypothesis (\ref{null-mean}) is in the space $(L^2({\cal I}), \, \| \cdot \|)$, i.e., $\mu_k=\mu_l$, for any pair of indices $(k,l)\in\{1,2,\ldots,K\}$, with $k \neq l$, means that  $ \|\mu_k-\mu_l\| =0$, and the alternative hypothesis (\ref{alt-mean}) means that 
$ \|\mu_k-\mu_l\| >0$.

\medskip
 
Such a testing problem can  be easily addressed   by  changing appropriately Step 3 of the bootstrap resampling algorithm of 
Section \ref{subsec.bootalg}. In particular, we replace   equation  
(\ref{eq.xboot}) in  Step 3 of this algorithm   by the following  equation
\begin{equation} \label{eq.xbootm}
X^+_{i,j}(t)=\overline{X}_{N}(t) + \epsilon^+_{i,j}(t), \quad  t \in {\cal I},
\end{equation}
where 
$$ \overline{X}_{N}(t) = \frac{1}{N}\sum_{i=1}^{K}\sum_{j=1}^{n_i} X_{i,j}(t), \quad t\in {\cal I},
$$ is the pooled mean function estimator and
$ \epsilon^+_{i,j}(t) =\widehat{\epsilon}_{i,J}(t)$, $t\in {\cal I}$,  where  $ J$ is a discrete random variable with $ P(J=j)=1/n_i$ for every  $ j =1,2, \ldots, n_i$, $i \in \{1,2,\ldots,K\}$. 
  Thus,  the bootstrap error functions $ \epsilon_{i,j}^+(t)$ for population $i$  appearing in equation (\ref{eq.xbootm}) are generated by randomly selecting a residual function from the set of estimated residual functions $ \widehat{\epsilon}_{i,j}(t)$ belonging to the same population $i \in \{1,2,\ldots,K\}$. This ensures that the covariance structure  of the functional observations in each population is retained by the bootstrap resampling  algorithm, which may be different for different populations, despite the fact that the bootstrap procedure  generates  $K$  populations of independent bootstrap  functional pseudo-observations  that have  the same mean function.
 In particular, conditional on the observed functional data ${\bf X}_N$, we have for the bootstrap functional pseudo-observations  $ X_{i,j}^+(t)$,  that, 
\[ \EE[X_{i,j}^+(t)]=\overline{X}_N(t),  \ \  t \in {\cal I},\]
and  
\begin{align*} 
Cov[X^+_{i,j}(t), X^+_{i,j}(s)]  &=   \EE[\epsilon_{i,j}^+(t)\epsilon_{i,j}^+(s)]
= \frac{1}{n_i} \sum_{j=1}^{n_i}\widehat{\epsilon}_{i,j}(t)\widehat{\epsilon}_{i,j}(s)\\
& =  \frac{1}{n_i} \sum_{j=1}^{n_i}(X_{i,j}(t)-\overline{X}_{i,n_i}(t)) (X_{i,j}(s)-\overline{X}_{i,n_i}(s)) \\
   & = \widehat{C}_{i,n_i}(t,s), \quad t,s \in {\cal I}.
\end{align*}

\begin{remark}  \label{re.equalmeanandcov}
 {\rm  
  Notice that if we use  $X^\circ_{i,j}(t)=\overline{X}_{N}(t) + \epsilon^*_{i,j}(t)$,  $ t \in {\cal I}$,  to generate the bootstrap pseudo-observations  instead of equation (\ref{eq.xbootm})
  with $\varepsilon^\ast_{i,j}(t)$ defined as in Step 3 of the algorithm in Section~\ref{subsec.bootalg}, then the   $X^\circ_{i,j}(t) $ will have in the $K$ groups an  identical mean function equal to $\overline{X}_N(t) $ {\it and} an identical covariance function equal to $ \widehat{C}_N(t,s)$. This may be of particular interest if one is interested in testing  {\em simultaneously}  the equality of mean functions and covariance functions between the $K$ populations.
  }
\end{remark}

 \begin{remark}  \label{re.adaptGauss}
 {\rm 
If  distributional assumptions  have been  imposed on the observed random functions $ X_{i,j}(t)$, then efficiency considerations  may suggest that such  assumptions 
 should be also taken into account in the implementation of the bootstrap resampling algorithms
 which are used to generated the bootstrap functional pseudo-observations. 
 For instance,  the assumption of 
 Gaussianity of the random paths $ X_{i,j}(t)$, $t\in {\cal I}$,  can  be incorporated in our bootstrap testing algorithm by allowing for the functional bootstrap pseudo-observations to follow a Gaussian processes on ${\cal I}$ with a mean and covariance function specified according the null hypothesis of interest.  
 }
\end{remark}

\section{{\sc  Bootstrap Validity}}
\label{sec:theor}

In this section,  we establish the validity of the introduced  bootstrap-based functional testing methodology applied to some test statistics recently proposed  in the literature for the important problems of testing the equality of  mean functions or  covariance functions between two populations.
 
\subsection{Testing the equality of two covariance functions}
\subsubsection{Test Statistics and Limiting Distributions}
\label{subsubsec:limited-T}

For testing the equality of two covariance operators, it is natural to evaluate the Hilbert-Schmidt norm of the difference of the corresponding sample covariance operators $\widehat{\mathcal C}_1$ and $\widehat{\mathcal C}_2$, defined as
$$
\widehat{\mathcal C}_i =\frac{1}{n_i}\sum_{j=1}^{n_i}(X_{i,j}-\overline{X}_ {i,n_i})\otimes (X_{i,j}-\overline{X}_ {i,n_i}), \quad i=1,2,
$$
or, equivalently, as 
$$
\widehat{\mathcal C}_i(f)(t) = \frac{1}{n_i} \sum_{j=1}^{n_i} \langle X_{i,j}-\overline{X}_{i,n_i}, f \rangle (X_{i,j}(t)-\overline{X}_{i,n_i}(t)), \quad f \in L^2({\cal I}), \quad t\in {\cal I},\quad i=1,2.
$$
Such an approach  has been recently proposed  by Boente {\em et al}. (2014)  by considering the test statistic 
\[ T_N=N\|\widehat{\mathcal C}_1-
\widehat{\mathcal C}_2\|_S^2, \]
where $N=n_1+n_2$.   If $ n_1/N\rightarrow \theta_1\in (0,1)$, $ \EE\|X_{i,1}\|^4 <\infty$,  $ i \in \{1,2\} $,  and  the null hypothesis  $H_0$ given in (\ref{null-cov}) with $K=2$ 
is true, then,  as $ n_1, n_2 \rightarrow \infty$,  Boente {\em et al.} (2014) showed  that  $ T_N$
 converges weakly to $\sum_{l=1}^{\infty} \lambda_l Z^2_l$,  where $ Z_l$ are independently distributed standard Gaussian random variables and  $\lambda_l$ are the eigenvalues  of the pooled operator ${\mathcal  B}=\theta_1^{-1}{\mathcal  B}_1 + (1- \theta_1)^{-1}{\mathcal B}_2$, where $ {\mathcal B}_i$ is  the covariance operator
 of the limiting Gaussian random element  $U_i$ to which $ \sqrt{n_i}(\widehat{\mathcal C}_i-{\mathcal C}_i)$ converges weakly as $ n_i\rightarrow \infty$, $i=1,2$. Since the limiting distribution of $T_N$ depends on the unknown infinite  eigenvalues $ \lambda_l$, $l\geq 1$, 
 implementation of this asymptotic result for calculating critical values of the test is difficult.
With this in mind, Boente {\em et al.} (2014) have proposed a  bootstrap calibration procedure of the distribution of the test statistic $T_N$.
\medskip

Another, related, approach for testing the equality of two covariance functions, is to evaluate the distance  between the sample covariance functions $\widehat{C}_{1,n_1}(t,s)$ and $\widehat{C}_{2,n_2}(t,s)$,  $t,s \in {\cal I}$, of  each group  and the pooled sample covariance function $ \widehat{C}_N(t,s)$, $t,s \in {\cal I}$, based on the entire set of functional observations. Therefore, looking at projections of  $\widehat{C}_{1,n_1}(t,s)-\widehat{C}_{2, n_2}(t,s)$, $t,s \in {\cal I}$, on certain directions reduces the dimensionality of the problem.   Such  approaches  have  been considered by Panaretos {\em et al.} (2010) (for Gaussian curves) and Fremdt {\em et al} (2012), (for non-Gaussian curves), where the asymptotic distributions of the corresponding test statistics proposed under the null hypothesis have been derived 

\medskip
More specifically,  denote by $(\widehat{\lambda}_k, \widehat{\varphi}_k)$, $k=1,2,\ldots,N$, the eigenvalues/eigenfunctions of the pooled sample covariance operator $\widehat{\cal C}_N$ defined by the kernel $\widehat{C}_N(t,s)$, i.e.,
\[  
\widehat{\lambda_k} \widehat{\varphi}_k(t) = \widehat{\cal C}_N(\widehat{\varphi}_k)(t)=
\int_0^1 \widehat{C}_N(t,s)\widehat{\varphi}_k(s)ds, \ \ \ t\in {\cal I}, \ \ \ k=1,2,\ldots,N,
\]
with $\widehat{\lambda}_1 \geq \widehat{\lambda}_2 \geq  \cdots$. 
(We can and will assume that the $\widehat{\varphi}_k(t)$, $k=1,2,\ldots, N$, $t \in {\cal I}$, form an orthonormal system.) 
Select  a natural number $p$ and consider,  for $ i=1,2, \ldots, p$, the projections  
\[ \widehat{a}_{k,j}(i) = <X_{k,j}-\overline{X}_{k,n_k}, \widehat{\varphi}_{i}> = \int_{\cal I} (X_{k,j}(t)-\overline{X}_{k,n_k}(t))\widehat{\varphi}_i(t)dt, \quad j=1,2, \ldots, n_k, \;k=1,2.
\]
For $ 1 \leq r,m \leq p$,  consider  the matrices $\widehat{A}_{k,n_k}$, $k=1,2$, with elements
\[ 
\widehat{A}_{k,n_k}(r,m) = \frac{1}{n_k}\sum_{j=1}^{n_k}\widehat{a}_{k,j}( r)\widehat{a}_{k,j}( m), \quad k=1,2.
\]
Notice that, for $ 1 \leq r,m \leq p$,  $\widehat{\Delta}_N(r,m) := \widehat{A}_{1,n_1}(r,m)-  \widehat{A}_{2,n_2}(r,m)$
is the projection of the difference $ \widehat{C}_{1,n_1}(t,s)-\widehat{C}_{2,n_2}(t,s)$ in the direction of $ \widehat{\varphi}_r(t) \widehat{\varphi}_m(s)$, $t, s \in {\cal I}$.

\medskip
Panaretos {\em et al.} (2010) considered then the  test statistic 
\begin{align*} 
T_{p,N}^{(G)}  & = \frac{n_1 n_2}{N} \sum_{1\leq r,m\leq p}\frac{\widehat{\Delta}_N^2(r,m)}{2\,\widehat{\lambda}_r\widehat{\lambda}_m}.
 \end{align*}
They showed, that, under  the assumption that $  X_{i,1}(t)$, $i \in \{1,2\}$, $ t \in {\cal I}$, are Gaussian  processes, and if  $ n_1,n_2\rightarrow \infty$ such that  $ n_1/N\rightarrow \theta_1\in (0,1)$, $ \EE\|X_{i,1}\|^4 <\infty$, $ i \in \{1,2\} $ and the null hypothesis $H_0$ given in (\ref{null-cov}) with $K=2$ is true, then, 
 $T_{p,N}^{(G)}$ converges weakly to a $\chi^2_{p(p+1)/2}$ distribution.

\medskip

The non-Gaussian  case has been recently investigated by Fremdt {\em et al.} (2012).  In particular,
they considered the matrix $\widehat{\Delta}_N=(\widehat{\Delta}_N(r,m))_{r,m=1,2, \ldots, p}$
 and defined 
$\widehat{\xi}_N= {\rm vech}(\widehat{\Delta}_N)$,
i.e.,  the vector containing the elements on and below the main diagonal of  $ \widehat{\Delta}_N$.
Fremdt {\em et al}. (2012) proposed then  the  test statistic  
\[ 
T_{p,N} = \frac{n_1 n_2}{N}\,  \widehat{\xi}_N^{T}\, \widehat{L}^{-1}_N\, \widehat{\xi}_N,
\]
where $ \widehat{L}_N$ is an estimator of the (asymptotic) covariance matrix of $ \widehat{\xi}_N$.
They showed that if $ n_1,n_2 \rightarrow \infty$ such that $ n_1/N\rightarrow \theta_1\in (0,1)$, $ \EE\|X_{i1}\|^4 <\infty$,  $ i \in \{1,2\} $,  and  the null hypothesis $ H_0$ given in (\ref{null-cov}) with $K=2$ is true, then,  $T_{2,N}$ converges weakly to  a $\chi^2_{p(p+1)/2}$ distribution.  Furthermore,   under the same set of assumptions,  consistency of the test  $T_{p,N}$  has been established under the alternative, that is when  the  covariance functions $C_1$ and $ C_2$ differ. 

\subsubsection{Consistency of the Bootstrap}

We  apply the bootstrap procedure introduced in Section~\ref{subsec.bootalg}  to approximate the distributions of the test statistics $T_N$, $T_{p,N}^{(G)}$ and of $ T_{p,N}$ under the null hypothesis. 
To this end, let $ X^\ast_{i,j}(t)$, $i \in \{1,2\}$, $j=1,2, \ldots, n_i$, $t \in {\cal I}$, be the bootstrap functional pseudo-observations generated according to this 
bootstrap procedure.
Let \[ T_N^\star=N\|\widehat{\mathcal C}_1^\star-
\widehat{\mathcal C}_2^\star\|_S^2, \]
where $\widehat{\mathcal C}_1^\star$ and $\widehat{\mathcal C}_2^\star$ are the sample covariance operators of the two groups but calculated using the the bootstrap functional pseudo-observations 
$ X^\ast_{i,j}(t)$, $i \in \{1,2\}$, $j=1,2, \ldots, n_i$, $t \in {\cal I}$.
Let 
\begin{align*} 
T_{p,N}^{\ast (G)}  & = \frac{n_1 n_2}{N} \sum_{1\leq r,m\leq p}\frac{1}{2}\frac{\widehat{\Delta}^{\ast^2}_N(r,m)}{
\widehat{\lambda}^\ast_r\widehat{\lambda}^\ast_m},
 \end{align*}
where $\widehat{\Delta}^{\ast^2}_N(r,m) $ and $ \widehat{\lambda}^\ast_r$ 
are the same statistics as $\widehat{\Delta}^{2}_N(r,m) $ and $ \widehat{\lambda}_r$ appearing in $ T_{p,N}^{(G)}$ but calculated using the bootstrap functional pseudo-observations 
$ X^\ast_{i,j}(t)$, $i \in \{1,2\}$, $j=1,2, \ldots, n_i$, $t \in {\cal I}$. Similarly, let
 \[ T_{p,N}^\ast = \frac{n_1 n_2}{N}\,  \widehat{\xi}_N^{\ast^T}\, \widehat{L}_N^{\ast^{-1}}\, \widehat{\xi}^\ast_N,
 \]
where $ \widehat{\xi}_N^{\ast} $ and $ \widehat{L}_N^{\ast} $ are the same statistics as  $ \widehat{\xi}_N $ and $ \widehat{L}_N  $ appearing in  $ T_{p,N}$ but calculated using the bootstrap functional pseudo-observations 
$ X^\ast_{i,j}(t)$, $i \in \{1,2\}$, $j=1,2, \ldots, n_i$, $t \in {\cal I}$. 
The following results are then true.

\begin{theorem}
\label{th:size-TNstar}
If  $ E\|X_{i,1}\|^4 <\infty$, $ i \in \{1,2\} $, and $ n_1/N\rightarrow \theta_1\in (0,1)$, then, as $ n_1, n_2 \rightarrow \infty$,
\[ 
\sup_{x\in \Re}\Big|\mathbb{P}\big(  T_{N}^\ast  \leq x  \mid {\bf X}_N\big) - \mathbb{P}_{H_0}\big(  T_{N}  \leq x  \big)\Big| \ \rightarrow \ 0, \quad \text{in probability,}
\] 
where $ \mathbb{P}_{H_0}(T_N\leq x)$, $ x \in \Re$, denotes the distribution function of $T_N$ 
when  $H_0$ given in (\ref{null-cov}) with $K=2$ is true and ${\mathcal  B}_1={\mathcal  B}_2$.
\end{theorem}

Notice that by the above theorem, the suggested bootstrap procedure  leads to consistent estimation of the critical values of the test $T_N$  for which the asymptotic approximation discussed in Section~\ref{subsubsec:limited-T}  is difficult to implement in practice.

\begin{theorem}
\label{th:size-T1star}
Assume that $  X_{i,1}(t)$, $i \in \{1,2\}$, $ t \in {\cal I}$, are Gaussian  processes. If  $ \EE\|X_{i,1}\|^4 <\infty$, $ i \in \{1,2\} $, and $ n_1/N\rightarrow \theta_1\in (0,1)$, then, as $ n_1, n_2 \rightarrow \infty$,
\[ 
\sup_{x\in \mathbb{R}}\Big|\mathbb{P}\big(  T_{p,N}^{\ast (G)}  \leq x  \mid {\bf X}_N\big) - \mathbb{P}_{H_0}\big(  T_{p,N}^{(G)}  \leq x  \big)\Big| \ \rightarrow \ 0, \quad \text{in probability,}
\]
where $ \mathbb{P}_{H_0}(T_{p,N}^{(G)}\leq x)$, $ x \in \Re$, denotes the distribution function of $T_{p,N}^{(G)}$ 
when  $H_0$ given in (\ref{null-cov}) with $K=2$ is true.
\end{theorem}

\begin{theorem}
\label{th:size-T2star}
If  $ \EE\|X_{i,1}\|^4 <\infty$, $ i \in \{1,2\} $, and $ n_1/N\rightarrow \theta_1\in (0,1)$, then, as $ n_1, n_2 \rightarrow \infty$,
\[ 
\sup_{x\in \mathbb{R}}\Big|\mathbb{P}\big(  T_{p,N}^\ast  \leq x  \mid {\bf X}_N\big) - \mathbb{P}_{H_0}\big(  T_{p,N}  \leq x  \big)\Big| \ \rightarrow \ 0, \quad \text{in probability,}
\] 
where $ \mathbb{P}_{H_0}(T_{p,N} \leq x)$, $ x \in \Re$, denotes the distribution function of $T_{p,N}$ 
when  $H_0$ given in (\ref{null-cov}) with $K=2$ is true.
\end{theorem}

\begin{remark}  \label{re.powcov}
 {\rm 
 Notice that if $H_1$ is true and $\|{\mathcal C}_1-{\mathcal C}_2\|_S>0$, then, we have that, as $n_1, n_2 \rightarrow \infty$, $ T_{N} \rightarrow \infty$, in probability. Theorem~\ref{th:size-TNstar} implies then that the test $T_{N}$ based on the bootstrap critical values obtained using the distribution of the test $ T_{N}^\ast$ is consistent, i.e., its power approaches unity,  as $ n_1,n_2 \rightarrow \infty$.
 Furthermore, if $H_1$ is true and if $ {\bf \xi} ={\rm vech}(D) \neq 0$, where $ D$ is the $p\times p$  matrix 
 $D =\big((\int_{{\mathcal I}}\int_{{\mathcal I}}(C_1(t,s)-C_2(t,s))\varphi_i(t)\varphi_j(s)dtds, 1\leq i,j\leq p\big)$, then, under the same assumptions as in Theorem 3 of Fremdt {\em et al.} (2012), we have that, as $n_1, n_2 \rightarrow \infty$, $ T_{p,N} \rightarrow \infty$, in probability. Theorem~\ref{th:size-T2star} implies then that, under these assumptions,  the test $T_{p,N}$ based on the bootstrap critical values obtained using the distribution of the test $ T_{p,N}^\ast$ is also consistent, i.e., its power approaches unity, as $ n_1,n_2 \rightarrow \infty$. 
 }
\end{remark}

\subsection{Testing the equality of two mean functions}
\subsubsection{Test Statistics and Limiting Distributions}
\label{subsubsec:limited-S}

For testing the equality of two mean functions, it is natural to compute the $L^2({\cal I})$-distance between the two sample mean functions $\overline{X}_{1,n_1}(t)$ and $\overline{X}_{2,n_2}(t)$, $t \in {\cal I}$. Such an approach was considered by Benko {\em et al.} (2009) and Horv\'ath \& Kokoszka (2012,  Chapter 5) using  the test statistic
$$
S_N=\frac{n_1 n_2}{N} \| \overline{X}_{1,n_1}-\overline{X}_{2,n_2} \|^2.
$$
If $ n_1/N\rightarrow \theta_1\in (0,1)$, $ \EE\|X_{i,1}\|^4 <\infty$,  $ i \in \{1,2\} $,  and  the null hypothesis $ H_0$ given in (\ref{null-mean}) with $K=2$ is true, then, as $ n_1, n_2 \rightarrow \infty$, they showed that $S_{N}$ converges weakly to $\int_{\cal I} \Gamma^2(t)dt$,
where $\{\Gamma(t):\, t \in {\cal I}\}$ is a Gaussian process satisfying 
$\EE[\Gamma(t)]=0$ and $\EE[\Gamma(t)\Gamma(s)]=(1-\theta)\, C_{1}(t,s)+\theta\, C_{2}(t,s)$, $t,s \in {\cal I}$.
They have also showed consistency of the test, in the sense that if the alternative hypothesis $ H_1$ given in (\ref{alt-mean})  is true, then, as $ n_1, n_2 \rightarrow \infty$, $S_N \rightarrow \infty$, in probability.

\medskip

Notice that the limiting distribution of the test statistic $S_N$ depends on the unknown covariance functions $C_{1}$ and $C_{2}$. Hence,  analytical 
calculation of critical values of this test turns out to be difficult in practice.  To overcome this problem, Horv\'ath \& Kokoszka (2012,  Chapter 5), considered two projections versions of the test statistic $S_N$. Note that, using the KLE, it follows that 
$
\Gamma(t) = \sum_{k=1}^{\infty} \sqrt{\tau_k}\, N_k\,\phi_k(t)$ and $\int_{\cal I} \Gamma^2(t)dt = \sum_{k=1}^{\infty} \tau_k \,N_k^2$, $t \in {\cal I}$,
where $N_k$, $k=1,2,\ldots$, is a sequence of independent standard Gaussian random variables, and $\tau_1 \geq \tau_2 \geq \ldots$ and $\phi_1(t),\phi_2(t),\ldots$, $t \in {\cal I}$, are the eigenvalues and eigenfunctions of the operator ${\cal Z}$ determined by the kernel $Z(t,s)=(1-\theta)\, C_{1}(t,s)+\theta\, C_{2}(t,s)$, $t,s \in {\cal I}$, $\theta \in (0,1)$. In view of this, Horv\'ath \& Kokoszka (2012,  Chapter 5) considered projections on the space determined by the $p$ leading eigenfunctions of the operator ${\cal Z}$. Assume that the eigenvalues of the operator ${\cal Z}$ satisfy $\tau_1 > \tau_2 >\ldots>\tau_p>\tau_{p+1}$, i.e., that there exists at least $p$ distinct (positive) eigenvalues of the operator ${\cal Z}$. Let 
$$
\hat{a}_i = <\overline{X}_{1,n_1}-\overline{X}_{2,n_2}, \hat{\phi}_i>=\int_{\cal I} (\overline{X}_{1,n_1}(t)-\overline{X}_{2,n_2}(t)) \hat{\phi}_i(t)dt, \quad i=1,2,\ldots,p,
$$
be the projection of the difference $\overline{X}_{1,n_1}(t)-\overline{X}_{2n_2}(t)$, $t \in {\cal I}$, into the linear space spanned by $\hat{\phi}_1(t)$, $\hat{\phi}_2(t)$,\ldots,$\hat{\phi}_p(t)$, $t \in {\cal I}$, the eigenfunctions related to the sample estimator $\hat{Z}_N(t,s)$ of the kernel $Z(t,s)$, $t, s \in {\cal I}$. Based on the above, Horv\'ath \& Kokoszka (2012,  Chapter 5) considered the following test statistics 
$$
S_{p,N}^{(1)}=\frac{n_1 n_2}{N} \sum_{k=1}^p \frac{\hat{a}_k^2}{{\hat \tau}_k} \quad \text{and} \quad
S_{p,N}^{(2)}=\frac{n_1 n_2}{N} \sum_{k=1}^p \hat{a}_k^2.
$$
If $ n_1/N\rightarrow \theta_1\in (0,1)$, $ \EE\|X_{i,1}\|^4 <\infty$,  $ i \in \{1,2\} $,  and  the null hypothesis $ H_0$ given in (\ref{null-mean}) with $K=2$ is true, then, as $ n_1, n_2 \rightarrow \infty$, they showed that, $S_{p,N}^{(1)}$ converges weakly to  a $\chi^2_{p}$-distribution while $S_{p,N}^{(2)}$ converges weakly to $\sum_{k=1}^p \tau_k N_k^2$.  Under the assumption that $\mu_1(t)-\mu_2(t)$, $t \in {\cal I}$, is not orthogonal to the linear span of $\phi_1(t),\phi_2(t),\ldots,\phi_p(t)$, $t \in {\cal I}$,  they have also showed consistency, in the sense that,  if $ n_1/N\rightarrow \theta_1\in (0,1)$, $ \EE\|X_{i,1}\|^4 <\infty$,  $ i \in \{1,2\} $, and the alternative hypothesis $H_1$ given in (\ref{alt-mean}) with $K=2$ is true, then, as $ n_1, n_2 \rightarrow \infty$, $S_{p,N}^{(1)} \rightarrow \infty$ and $S_{p,N}^{(2)} \rightarrow \infty$, in probability.

\subsubsection{Consistency of the Bootstrap}
To approximate the distribution of the test statistics  $S_N$,  $S_{p,N}^{(1)}$ and $S_{p,N}^{(2)}$, we apply the bootstrap procedure proposed in Section~\ref{re.adapt}. 
To this end, let $ X^+_{i,j}(t)$, $i \in \{1,2\}$, $j=1,2, \ldots, n_i$, $t \in {\cal I}$, be the bootstrap functional pseudo-observations generated according to this bootstrap algorithm 
and define
\begin{align*} 
S_N^+ &=\frac{n_1 n_2}{N} \| \overline{X}^+_{1,n_1}-\overline{X}^+_{2,n_2} \|^2,\\
S_{p,N}^{+ (1)}  & = \frac{n_1 n_2}{N} \sum_{k=1}^p \frac{\hat{a}_k^{+^2}}{{\hat \tau}_k^{+}}, \;\; \text{and}\\
S_{p,N}^{+ (2)}  & =\frac{n_1 n_2}{N} \sum_{k=1}^p {\hat a}_k^{+^2}
\end{align*}
 be the same statistic as $S_N$, $S_{p,N}^{(1)}$ and $S_{p,N}^{(2)}$, respectively, but calculated using the bootstrap functional pseudo-observations $ X^+_{i,j}(t)$, $i \in \{1,2\}$, $j=1,2, \ldots, n_i$, $t \in {\cal I}$.
 The following results are then true.

\begin{theorem}
\label{th:size-Ustar}
If  $ \EE\|X_{i,1}\|^4 <\infty$, $ i \in \{1,2\} $, and $ n_1/N\rightarrow \theta_1\in (0,1)$, then, as $ n_1, n_2 \rightarrow \infty$,
\[ 
\sup_{x\in \mathbb{R}}\Big|\mathbb{P}\big(  S_{N}^+  \leq x  \mid {\bf X}_N\big) - \mathbb{P}_{H_0}\big(  S_{N}  \leq x  \big)\Big| \ \rightarrow \ 0, \quad \text{in probability,}
\] 
where $ \mathbb{P}_{H_0}(S_N \leq x)$, $ x \in \Re$, denotes the distribution function of 
$S_N$ when  $H_0$ given in (\ref{null-mean}) with $K=2$ is true.
\end{theorem}

Notice that by  the above theorem, the suggested bootstrap procedure  leads to consistent estimation of the critical values of the test $S_N$  for which the asymptotic approximations discussed in Section~\ref{subsubsec:limited-S}  is difficult to implement in practice.

\begin{theorem}
\label{th:size-S1star}
If  $ \EE\|X_{i,1}\|^4 <\infty$, $ i \in \{1,2\} $, and $ n_1/N\rightarrow \theta_1\in (0,1)$, then, as $ n_1, n_2 \rightarrow \infty$,
$$
(i)\quad \quad \sup_{x\in \mathbb{R}}\Big|\mathbb{P}\big(  S_{p,N}^{+ (1)}  \leq x  \mid  {\bf X}_N\big) - \mathbb{P}_{H_0}\big(  S_{p,N}^{(1)}  \leq x  \big)\Big| \ \rightarrow \ 0, \quad \text{in probability},
$$
and
$$ 
(ii)\quad \quad\sup_{x\in \mathbb{R}}\Big|\mathbb{P}\big(  S_{p,N}^{+ (2)}  \leq x  \mid {\bf X}_N\big) - \mathbb{P}\big(  S_{p,N}^{(2)}  \leq x  \big)\Big| \ \rightarrow \ 0, \quad \text{in probability},
$$
where $ \mathbb{P}_{H_0}(S_{p,N}^{(1)} \leq x)$ and $ \mathbb{P}_{H_0}(S_{p,N}^{(2)} \leq x)$, $ x \in \Re$, denote the distribution functions of $S_{p,N}^{(1)}$ and $S_{p,N}^{(2)}$, respectively, when  $H_0$ given in (\ref{null-mean}) with $K=2$ is true.
\end{theorem}

\begin{remark}  \label{re.powmean}
 {\rm 
 Under the same assumptions as in Theorem 5.2 of   Horv\'ath \& Kokoszka (2012) and if $ \|\mu_1-\mu_2\| >0$ we have that  $S_N \rightarrow \infty$ in probability. This result together with  Theorem~\ref{th:size-Ustar} imply consistency of the test $S_N$ using the bootstrap critical values of the distribution of the test $S^+_N$, i.e., its power approaches unity, as $ n_1,n_2 \rightarrow \infty$. Furthermore, if the difference  $ \mu_1-\mu_2$ is not orthogonal to the linear span of $ \varphi_1, \varphi_2, \ldots, \varphi_p$, then, as $n_1, n_2 \rightarrow \infty$, $ S_{p,N}^{(1)} \rightarrow \infty$ and $S_{p,N}^{(2)} \rightarrow \infty$, in probability.  Theorem~\ref{th:size-S1star} implies then that, under these conditions, the tests $S_{p,N}^{(1)}$ and $S_{p,N}^{(2)}$ based on the bootstrap critical values of the distributions of the tests $ S_{p,N}^{+ (1)}$ and $ S_{p,N}^{+ (2)}$, respectively, are also consistent, i.e., their power approaches unity, as $ n_1,n_2 \rightarrow \infty$..  
 }
\end{remark}

\section{{\sc Numerical Results}}
\label{sec:num}
In this section, we  evaluate the finite sample behavior of the proposed  bootstrap-based functional testing procedures, for testing the equality of two covariance functions or the equality of two mean functions, by means of several  simulations and  compare our results with those  based on  classical asymptotic approximations of the distribution of the test statistic. An illustration 
to an interesting real-life dataset is also presented.

\subsection{\sc Simulations }
\label{subsec:s-s}
Following  Fremdt {\em et al.} (2012), we have simulated Gaussian curves $X_1(t)$ and $X_2(t)$, $t \in {\cal I}$,  as Brownian motions (BM) or Brownian bridges (BB), and non-Gaussian (NG) curves $X_1(t)$ and $X_2(t)$, $t \in {\cal I}$, via
\begin{equation}
X_i(t)=A \sin(\pi t) + B \sin(2\pi t) + C \sin(4 \pi t), \quad t \in {\cal I}, \ \  i\in \{1,2\},
\label{eq:NG-sim}
\end{equation}
where $A=7Y_1$, $B=3Y_2$, $C=Y_3$ with $Y_1$, $Y_2$ and $Y_3$ are  independent $t_5$-distributed random variables.  All curves were simulated at 500 equidistant points in the unit interval ${\cal I}$, and transformed into functional objects using the Fourier basis with 49 basis functions. For each data generating process, we considered 500 replications.  For practical and computational reasons, we have concentrated our analysis to sample sizes ranging from 
 $n_1=n_2=25$ to  $n_1=n_2=100$ random curves in each group, using also the three most common nominal levels $\alpha$, i.e., $\alpha \in \{0.01, 0.05, 0.10\}$. All bootstrap calculations are based on $B=1000$ bootstrap replications.

\medskip
We first illustrate the quality of the asymptotic and of the bootstrap approximations to the distribution of interest by considering  the test statistic $T_{2,N}$. For this, we first estimate the exact distribution of this test statistic under the null hypothesis by generating $10,000$ replications of functional data  $ {\bf X_N}$ using (\ref{eq:NG-sim}) and $n_1=n_2=25$ observations. We then compare the kernel density estimate of this exact distribution (obtained using a Gaussian kernel with bandwidth equal to 0.45) with that of the asymptotic $ \chi^2_3$ distribution and that of the bootstrap approximation using the algorithm described in Section~\ref{subsec.bootalg}  to generate the bootstrap pseudo-functional data $ {\bf X}_N^\ast$. Figure 6.1 and 6.2 presents the results obtained by applying the bootstrap to five randomly selected samples $ {\bf X}_N$. 

\begin{center}
{\bf Please insert Figure 6.1 and Figure 6.2 about here}
\end{center}

As it seen from these exhibits, which  present  density  estimates of the distribution of interest and corresponding QQ-plots, the density of the asymptotic $ \chi^2_3$ distribution does not provide an accurate approximation of the exact density of interest. In particular, it overestimates the exact density in  the crucial region of the right tail of this distribution, as it is clearly seen in Figure 6.1. This overestimation implies that the $T_{2,N}$ test using $\chi^2$
critical values will   lead to rejection rates that are below the desired nominal size of $ \alpha$, that is the test will be conservative.   On the other hand, and compared to the $\chi^2$ approximation, the bootstrap estimations  are much more  accurate and provide a very good approximation of  the exact distribution of the test statistic 
considered.   This behavior is also clearly seen in the QQ-plots presented in  Figure 6.2 .     

\medskip

We next investigate the  sizes behavior of the tests $T_{p,N}^{(G)}$ and $ T_{p,N}$  for testing equality of two covariance functions, both for BM and NG data, using the asymptotic $\chi^2$-approximation, where the corresponding tests are denoted  by $T_{p,N}^{(G)}$-Asym and $T_{p,N}$-Asym, respectively and their bootstrap approximations, where the corresponding tests are denoted by $T_{p,N}^{(G)}$-Boot and $T_{p,N}$-Boot, respectively. We have also tested the performance of the bootstrap approximation test $T_N^\star$, denoted by $T_N$-Boot. We use  either $n_1=n_2=25$,  $n_1=n_2=50$ or
$n_1=n_2=100$ curves, with either two ($p=2$) or three ($p=3$)  FPC's to perform the tests. 

\begin{center}
{\bf Please insert Table~\ref{tab:empsize-cov}   about here}
\end{center}

\begin{figure}
\begin{center}
\includegraphics[angle=0,height=8cm,width=15cm]{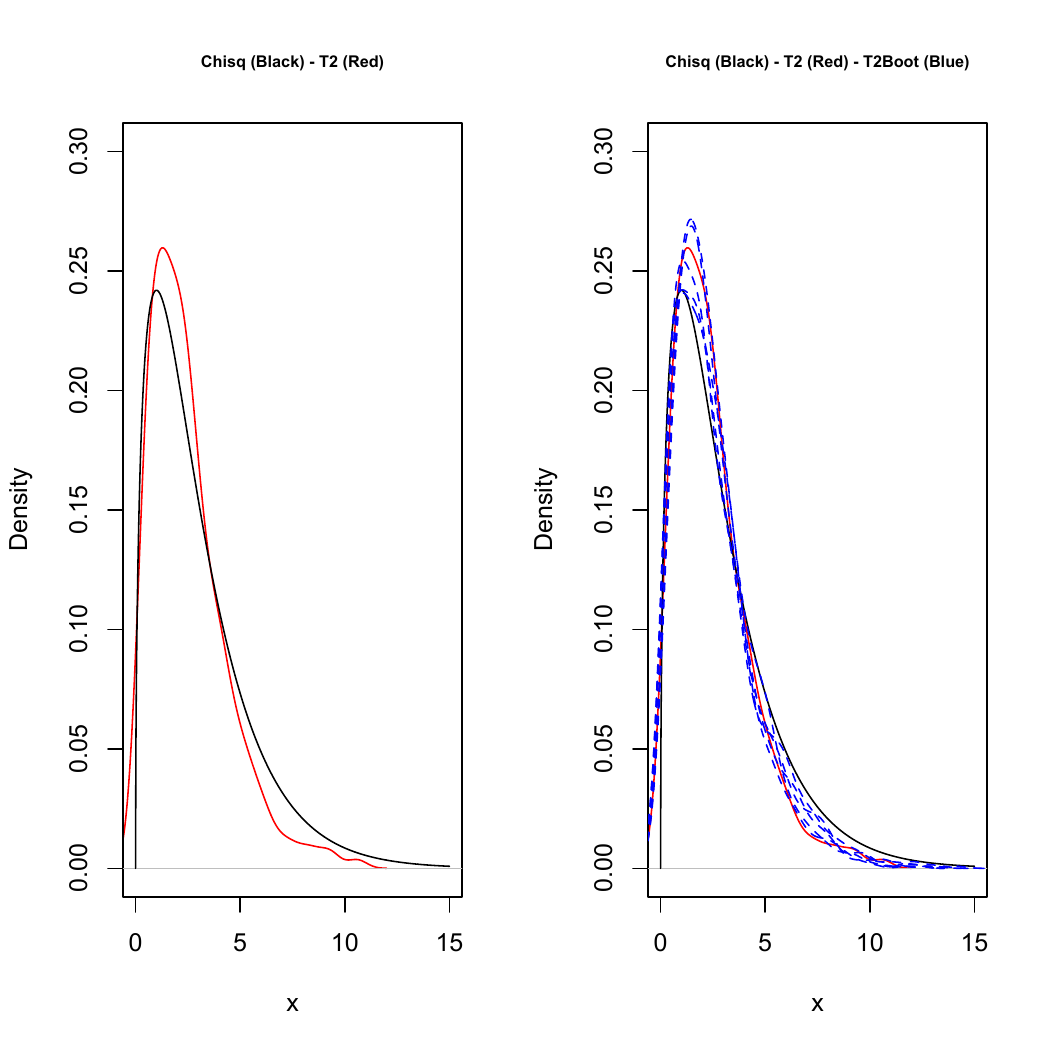}
\end{center}
\caption{Density estimates,  test statistic $T_{p,N}$-Asym: $n_1=n_2=25$, $p=2$.
Estimated exact density (red line),  $\chi^2_3$ approximation (black line) and 5  bootstrap approximations (dashed blue lines).}
\label{Fig:1}
\end{figure}

\begin{figure}
\begin{center}
\includegraphics[angle=0,height=10cm,width=15cm]{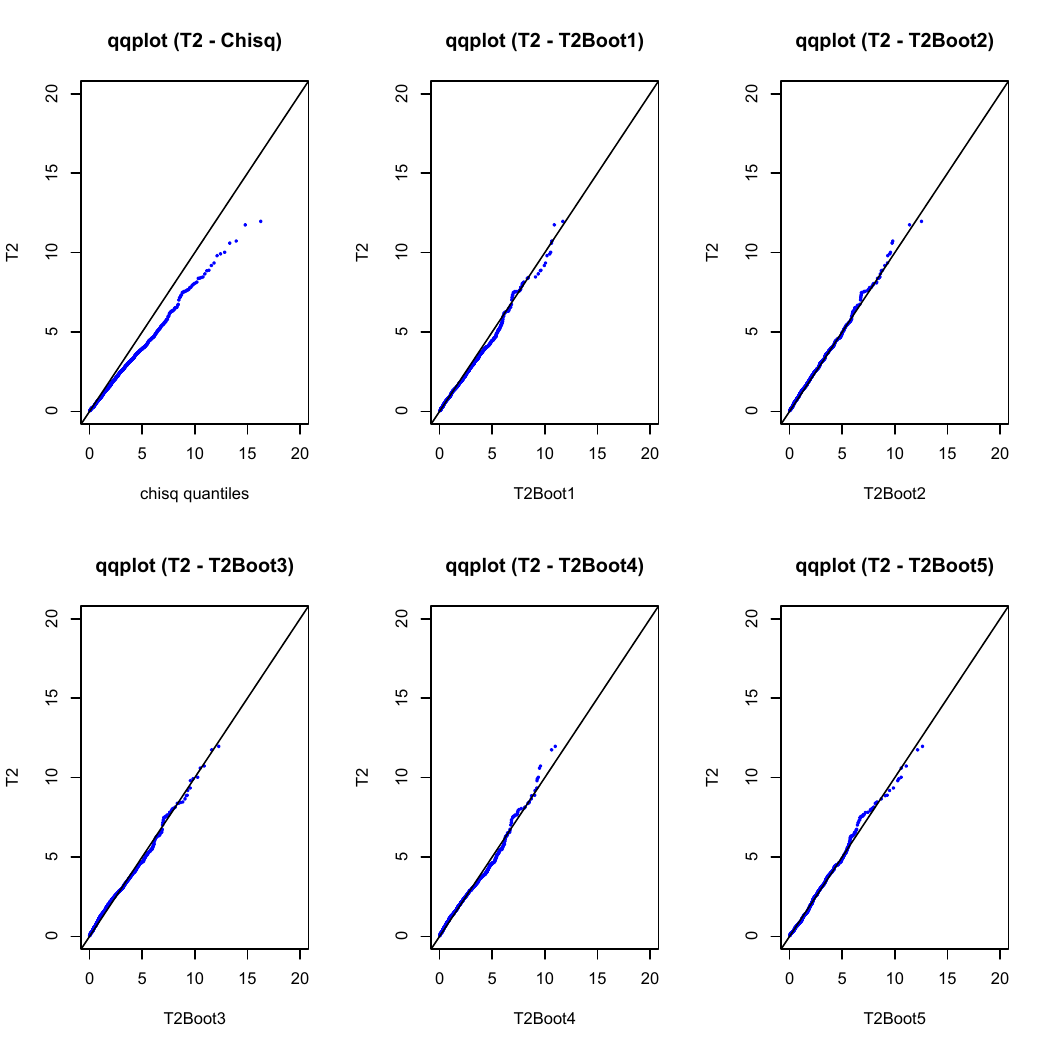}
\end{center}
\caption{QQ-plots,  Test statistic $T_{p,N}$-Asym: $n_1=n_2=25$, $p=2$.  Simulated exact distribution against the $ \chi^2_3$ distribution  (qqplot (T2-Chisq)) and against 5 bootstrap approximations (qqplot(T2-T2Bootj), $j=1,2,\ldots,5$.}
\label{Fig:2}
\end{figure}

Table~\ref{tab:empsize-cov} shows the empirical sizes obtained. As it   is evident from this table,   the test $T_{p,N}^{(G)}$-Asym has a severely inflated size in the case of NG data, due to violation of the assumption of normality,  a  behavior which was pointed out also in the simulation study of  Fremdt {\em et al.} (2012, Section 4). It is also evident that the test $T_{p,N}$-Asym has a severely under-estimated  size, confirming the visual evidence of Figure \ref{Fig:1} and Figure \ref{Fig:2}. 
On the other hand, it is clear that the tests $T_{p,N}^{(G)}$-Boot  and $ T_{p,N}$-Boot based on the bootstrap approximations have a very good  size behavior and do not suffer from the over- or  under- rejection problems from which the tests $T_{p,N}^{(G)}$ and $T_{p,N}$  based on asymptotic approximations suffer.  Notice that 
this very good size behavior is true even  for sample sizes as small as $n_1=n_2=25$ observations that have been considered in the simulation study. 
Finally, and in contrast to the behavior of the test $T_{p,N}^{(G)}$-Asym,  notice the  nice robustness property of the bootstrap-based counterpart  test $T_{p,N}^{(G)}$-Boot against deviations from Gaussianity in the case of  NG data. The advantage of the bootstrap can be seen here in the overall better performance of the bootstrap-based test $T_{N}$-Boot. Recall that the test $T_N$ does not requires  the 
choice of a truncation parameter $p$, referring to the number of  FPC's considered, and that  the 
asymptotic derivations  of the null distribution of this test  lead to results that are difficult to implement in order to calculate the critical values.
\medskip

To continue our  investigations of the finite sample behavior of the test statistics considered using  asymptotic and bootstrap approximations, we 
investigate the power properties  of the test $ T_{p,N}$ for the case of  NG data.  Due to the severe
 size distortions of the test $T_{p,N}^{(G)}$ for NG data, we do not include this  test statistic in our power study. 
We thus calculated the empirical rejection rates of  the tests $T_{p,N}$-Asym and $T_{p,N}$-Boot over $500$ replications, generated  for  either $n_1=n_2=25$ or $n_1=n_2=50$ observations and  two ($p=2$) FPC's.
  The curves in the first sample were generated according to (\ref{eq:NG-sim}) while the curves in the second sample were generated according to a scaled version of (\ref{eq:NG-sim}), i.e., $X_2(t)=\gamma X_1(t)$, $t \in {\cal I}$. The results are displayed for a selection of values of the scaling parameter $\gamma$, i.e., $\gamma \in \{2.0, 2.2, 2.4, 2.6, 2.8, 3.0\}$ in Table~\ref{tab:emppower-cov}. 
  
 \begin{center}
{\bf Please insert Table~\ref{tab:emppower-cov},  about here}
\end{center}

As it is   evident from Table~\ref{tab:emppower-cov},  the test $T_{p,N}$-Boot  based on bootstrap approximations 
 has  a much higher power that the test 
  $T_{p,N}$-Asym  based on asymptotic approximations. The low power of the test $T_{p,N}$-Asym is due to the  overestimation of the right-tail of the true density, as demonstrated in Figure 6.1. Notice that,  while this overestimation  
  leads to a conservative test under the null hypothesis, it leads to  a loss of  power  under the alternative.
 As can be expected, the power of the test $ T_{p,N}$-Boot improves  as the deviations from the null become larger 
 (i.e., larger values of  $\gamma$'s) and/or  as the  sample sizes increase.  Thus, and as our empirical evidence shows, the tests based on bootstrap approximations not only have  a better size behavior under the null hypothesis than those based on asymptotic approximations,  but they also have  a much  better power performance under the alternative.  However, the test $T_{N}$-Boot, which does not require the choice of $p$ (the truncation parameter), is overall the most powerful one. This clearly demonstrates the advantages of this bootstrap procedure.

\medskip
 
We next consider the finite sample  size and power  properties   of the asymptotic and of the bootstrap-based tests  considered for testing equality of mean functions.
 Table \ref{tab:empsize-mean} shows the empirical sizes of the tests for the equality of two mean functions, based on the statistics $S_{p,N}^{(1)}$-Asym and $S_{p,N}^{(2)}$-Asym (asymptotic approximations) and $S_{p,N}^{(1)}$-Boot, $S_{p,N}^{(2)}$-Boot and $S_{N}$-Boot (bootstrap approximations), both for BB and NG data.
 Sample sizes of either $n_1=n_2=25$ or $n_1=n_2=50$ have been considered with either two ($p=2$) or three ($p=3$)  FPC's and $500$ replications. 
 
  \begin{center}
{\bf Please insert Table~\ref{tab:empsize-mean}   about here}
\end{center}

As its evident from this table,  the tests $S_{p,N}^{(1)}$-Boot and $S_{p,N}^{(2)}$-Boot have sizes that are quite close to the nominal ones. The same is also for the tests $S_{p,N}^{(1)}$-Asym and $S_{p,N}^{(2)}$-Asym, although, in most of the cases, the empirical sizes of these tests exceed the nominal ones.
The advantage of the bootstrap can be seen here in the overall better performance of the bootstrap-based test $S_{N}$-Boot. Recall that the test $S_N$ does not requires  the 
choice of a truncation parameter $p$, referring to the number of  FPC's considered, and that  the 
asymptotic derivations  of the null distribution of this test  lead to results that are difficult to implement in order to calculate the critical values.

\medskip

Finally, we investigate  the power behavior of the tests considered.   For this, the empirical rejection frequencies  of the tests $S_{p,N}^{(1)}$-Asym, $S_{p,N}^{(2)}$-Asym, $S_{p,N}^{(1)}$-Boot, $S_{p,N}^{(2)}$-Boot and $S_{N}$-Boot have been calculated over $500$ replications  using  NG data, either $n_1=n_2=25$ or $n_1=n_2=50$ curves, with two ($p=2$) FPC's.   The curves in the two samples were generated according to the model $X_i(t)=\mu_i(t)+\epsilon_i(t)$ with $\epsilon_i(t)$ generated  according to (\ref{eq:NG-sim}), for  $i \in \{1,2\}$ and  $t \in {\cal I}$. The mean functions were set equal to $\mu_1(t)=0$ and $\mu_2(t)=\delta$ for each group respectively. The results obtained are displayed for a selection of values of the shift parameter $\delta$, i.e., $\delta \in \{1.0, 1.2, 1.4, 1.6, 1.8, 2.0\}$ in
 Table \ref{tab:emppower-mean}. 
 
   \begin{center}
{\bf Please insert Table~\ref{tab:emppower-mean}   about here}
\end{center}

As this table shows,  the  power results  for   the asymptotic-based tests $S_{p,N}^{(1)}$-Asym and $S_{p,N}^{(2)}$-Asym confirm  the findings of Horv\'ath \& Kokoszka  (2012, Table 5.1) which  have been obtained for larger sample sizes and for different deviations from the null.  Furthermore, 
the tests $S_{p,N}^{(1)}$-Boot and $S_{p,N}^{(2)}$-Boot show similar power behavior, although the slight better power performance of the asymptotic tests are due to the fact that these tests overestimate the nominal size, as mentioned above. The test $S_{N}$-Boot, which does not require the choice of $p$ (the truncation parameter), is the most powerful one. This clearly demonstrates the advantages of this bootstrap procedure.


\subsection{\sc Mediterranean Fruit Flies}
\label{subsec:r-d}
We now apply  the suggested bootstrap-based testing procedure  to a data set consisting of egg-laying trajectories of Mediterranean fruit flies ({\em Ceratitis capitata}), or medflies for short. This data set has been proved popular in the biological and statistical literature; see M\"uller \& Stadtm\"uller (2005) and references therein. It has also been analyzed by, e.g., Horv\'ath \& Kokoszka  (2012, Chapter 5) (for testing the equality of two mean functions) and by Fremdt {\em et al.} (2013) (for testing the equality of two covariance functions).

\medskip

We consider $N=534$ egg-laying curves of medflies who lived at least 43 days, but,  as in, e.g., Horv\'ath \& Kokoszka  (2012, Chapter 5) and Fremdt {\em et al.} (2013), we only consider the egg-laying activities on the first 30 days. Two versions of these egg-laying curves are considered and are scaled such that the corresponding curves in either version are defined on the interval ${\cal I}=[0,1]$. The curves in {\em Version 1} are denoted by $X_i(t)$ and represent the absolute counts of eggs laid by fly $i$ on day $\lfloor30t\rfloor$. The curves in {\em Version 2} are denoted by $Y_i(t)$ and represent the counts of eggs laid by fly $i$ on day $\lfloor30t\rfloor$ relative to the total number of eggs laid in the lifetime of  fly $i$. Furthermore, the 534 flies are classified into short-lived flies (those who died before the end of the 43rd day after birth) and long-lived flies (those who lived 44 days or longer). In this particular data set analyzed, there are $n_1=256$ short-lived flies and $n_2=278$ long-lived flies.  

\begin{figure}[h]
\begin{center}
\includegraphics[angle=0,height=6cm,width=7.5cm]{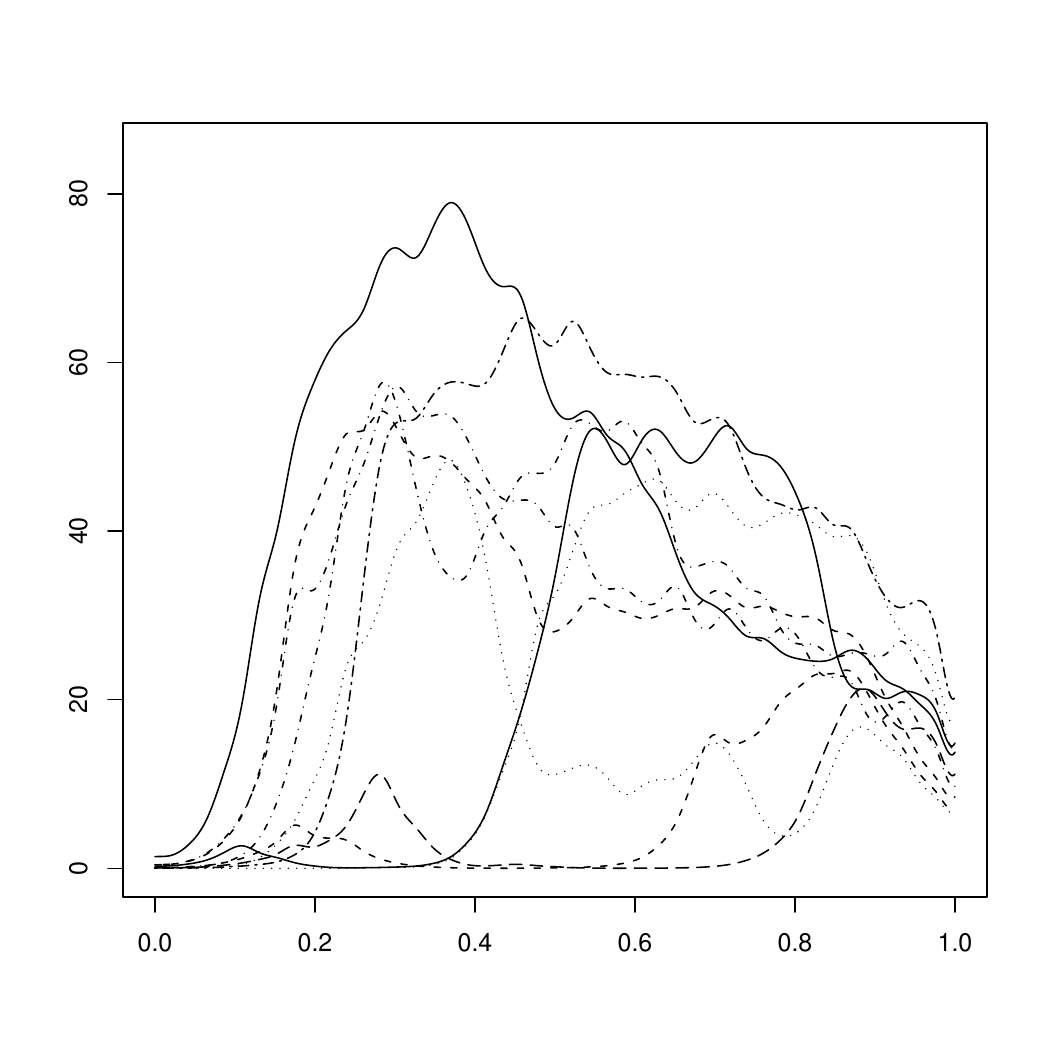}
\includegraphics[angle=0,height=6cm,width=7.5cm]{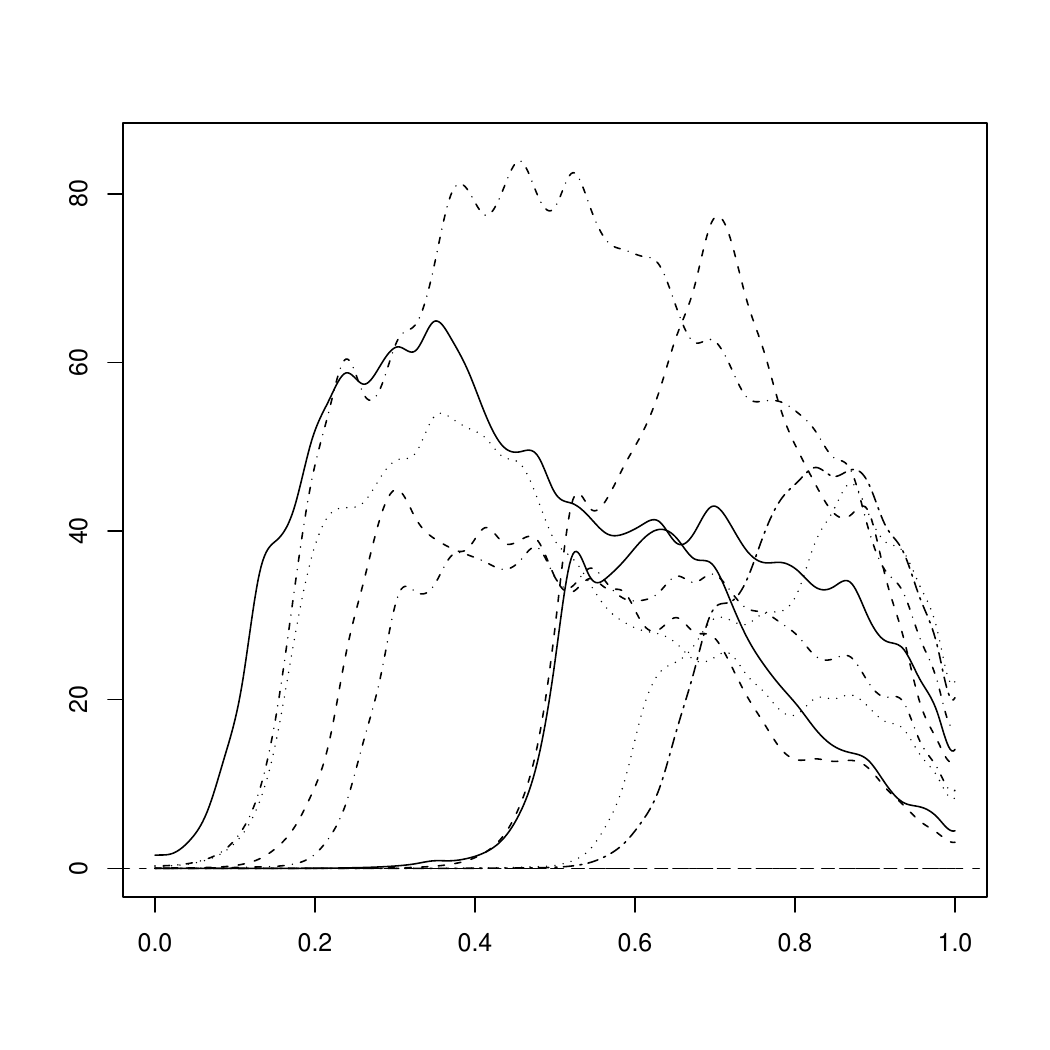}
\end{center}
\caption{10 randomly selected (smoothed) curves of short-lived (left panel) and 10 randomly selected (smoothed) curves of long-lived flies (right panel),  scaled on the interval ${\cal I}=[0,1]$.}
\label{fig:abs}
\end{figure}

\begin{figure}[h]
\begin{center}
\includegraphics[angle=0,height=6cm,width=7.5cm]{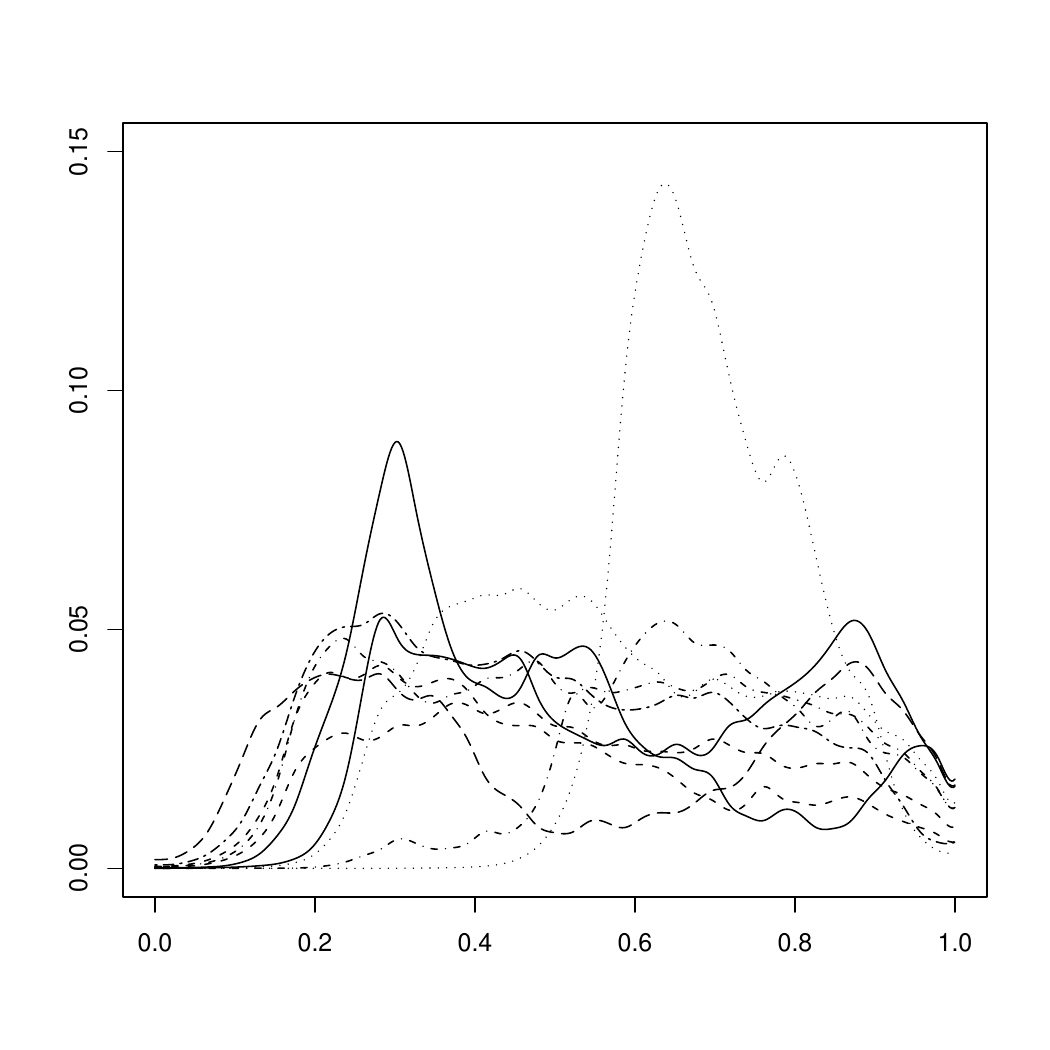}
\includegraphics[angle=0,height=6cm,width=7.5cm]{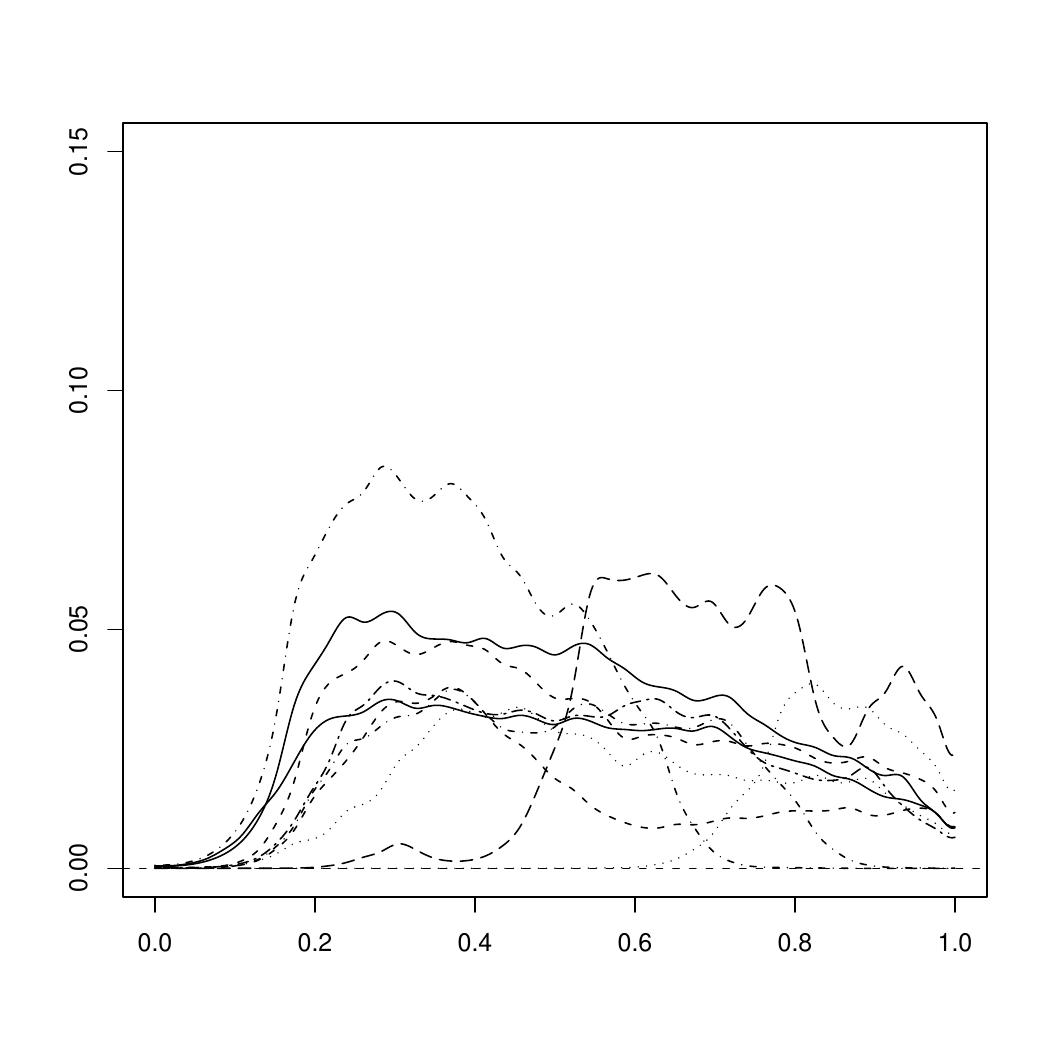}
\end{center}
\caption{10 randomly selected (smoothed) curves of short-lived (left panel) and 10 randomly selected (smoothed) curves of long-lived flies (right panel) - relative to the number of eggs laid in the fly's lifetime - scaled on the interval ${\cal I}=[0,1]$.}
\label{fig:rel}
\end{figure}

\medskip
Based on the above classification, we consider 2 samples. {\em Sample 1} represents the egg-laying curves of the {\em short-lived flies} $\{X_{1,i}(t):\; t \in {\cal I},\; i=1,2,\ldots, 256\}$ (absolute curves) or  $\{Y_{1,i}(t):\; t \in {\cal I},\; i=1,2,\ldots, 256\}$ (relative curves). {\em Sample 2} represents the egg-laying curves of the {\em long-lived flies} $\{X_{2,i}(t):\; t \in {\cal I},\; i=1,2,\ldots, 278\}$ (absolute curves) or  $\{Y_{2,i}(t):\; t \in {\cal I},\; i=1,2,\ldots, 278\}$ (relative curves). The actual curves were very irregular, hence originally smoothed slightly to produce the considered curves. Figure \ref{fig:abs} shows 10 randomly selected (smoothed) curves of short-lived and long-lived flies for {\em Version 1} while Figure \ref{fig:rel} shows 10 randomly selected (smoothed) curves of short-lived and long-lived flies for {\em Version 2}. The tests have been applied to such smooth curves, using a Fourier basis with with 49 basis functions for the representation into functional objects. Again, all bootstrap calculations  are based on $B=1000$ bootstrap replications.

\medskip

Table \ref{tab:medfly-cov} shows the $p$-values for the absolute (Figure \ref{fig:abs}) and the relative (Figure \ref{fig:rel}) egg-laying curves of the tests for the equality of covariance functions, using  the statistics $T_{1,N}$-Asym and $T_{2,N}$-Asym (based on asymptotic approximations) and $T_{1,N}$-Boot, $T_{2,N}$-Boot and $T_N$-Boot (based on bootstrap approximations).  

 \begin{center}
{\bf Please insert Table~\ref{tab:medfly-cov}   about here}
\end{center}

According to these results,   both  tests $T_{p,N}$-Asym and $T_{p,N}$-Boot show in the case of the absolute egg-laying curves a  uniform behavior across the range of the different  values of  $p$  that explain at least 85\% of the sample variance (a commonly used  rule-of-thumb recommendation). In particular,  and at the commonly used  $\alpha$-levels,
 the hypothesis of equality of the covariance functions cannot be rejected. However,  the opposite is true for the relative  egg-laying curves for which  the hypothesis of equality of the  covariance functions  should  be rejected   at the commonly used $\alpha$-levels and for most of the values of $p$ considered. Notice that the  bootstrap-based test $T_{p,N}$-Boot shows  in this case  a more stable behavior compared to  the test $T_{p,N}$-Asym, which does not  reject the null hypothesis for large values of $p$.    Furthermore, and due to the non-Gaussianity of the medfy data,  see Fremdt {\em et al.} (2013, Figure 3), and the over rejection of the test $T_{p,N}^{(G)}$-Asym demonstrated in Table ~\ref{tab:empsize-cov},  the results obtained using this  test for the absolute egg-laying  curves, provide little evidence  that the null hypothesis of equality of covariances should be rejected. The results using the test $T_{p,N}^{(G)}$-Boot are, however, more consistent with those obtained using the test $T_{p,N}$, designed for NG data.  
It is worth mentioning that the  bootstrap-based test $T_N$-Boot,  which does not require the choice of a truncation parameter $p$, leads to a clear rejection (for the relative egg-laying curves) and non-rejection for the absolute egg-laying curves) of the null hypothesis that the covariance functions are equal, demonstrating its usefulness in practical applications. 
 
\medskip
   \begin{center}
{\bf Please insert Table~\ref{tab:medfly-mean}   about here}
\end{center}

Table \ref{tab:medfly-mean} shows the $p$-values for the absolute (Figure \ref{fig:abs}) egg-laying curves of the tests for the equality of mean functions. The tests used are  $S_{p,N}^{(1)}$-Asym  and $S_{p,N}^{(2)}$-Asym (based on asymptotic approximations) and $S_{p,N}^{(1)}$-Boot, $S_{p,N}^{(2)}$-Boot and $S_N$-Boot (based on bootstrap  approximations). As it is evident from this table, the tests $S_{p,N}^{(2)}$-Asym and $S_{p,N}^{(2)}$-Boot provide a uniform behavior across the range of the first $p$ FPC's that explain at least 85\% of the sample variance, pinpointing to a rejection of the hypothesis  of equality of the mean functions. The tests  $S_{p,N}^{(1)}$-Asym and $S_{p,N}^{(1)}$-Boot show a more  erratic behavior, leading to rejection of the null hypothesis in the case of small or large $p$ and to non-rejection in the case of moderate $p$. This erratic behavior of the test $S_{p,N}^{(1)}$-Asym with respect to  the truncation parameter $p$ was also pointed out in Horv\'ath \& Kokoszka  (2012, Table 5.2). 
It is worth mentioning that the  bootstrap-based test $S_N$-Boot,  which does not require the choice of a truncation parameter $p$, leads to a clear rejection of the  null hypothesis that the mean functions are equal, demonstrating its usefulness in practical applications.  Since for the relative egg-laying curves the null hypothesis of equality of the two mean functions is rejected by all test statistics, we report only the results for the absolute egg-laying curves.

\section{\sc Conclusions} 
\label{sec:conc}
\medskip

We investigated properties of a simple bootstrap-based functional testing methodology which has been applied to the important problem of comparing the mean functions and/or the covariance functions  between several populations.  We theoretically justified the consistency of this bootstrap testing methodology applied to some tests statistics recently proposed in the literature,  and also demonstrated a very good size and power behavior in finite samples.\\

Although we restricted our theoretical investigations to some statistics recently proposed in the literature that build  upon the empirical  FPC's, the suggested bootstrap-based functional testing methodology can potentially be applied to other test statistics  too.  Such test statistics could be, for instance,  the likelihood ratio-type statistic for  testing the equality of two covariance functions considered in Gaines {\em et al.} (2011)  or  the regularized-based $M$-test considered in Kraus \& Panaretos (2012) for the same problem. 
Also, subject to appropriate modifications, we conjecture that the suggested basic resampling algorithm can be adapted  to different testing problems related to the comparisons of population characteristics like 
testing equality of distributions; see for instance the approaches of Hall \& Van Keilegom (2007) and Benko {\em et al.} (2009).
However, all the above investigations  require careful attention that is beyond the scope of the present work.

\section{{\sc Appendix: Proofs}}
We first fix some notation. Let $\widehat{\Delta}_N^\ast = \widehat{A}_{1,n_1}^\ast-\widehat{A}_{2,n_2}^\ast =\big( (\widehat{A}_{1,n_1}^\ast(r,m)-
\widehat{A}_{2,n_2}^\ast(r,m)), 1\leq r , m \leq p \big) $, where 
\[\widehat{A}_{1,n_1}^\ast(r,m)= \frac{1}{n_1}\sum_{j=1}^{n_1} \widehat{a}^\ast_{1,j}(r)\widehat{a}^\ast_{1,j}(m), \ \ \ \widehat{A}_{2,n_2}^\ast(r,m)= \frac{1}{n_2}\sum_{j=1}^{n_2} \widehat{a}^\ast_{2,j}(r)\widehat{a}^\ast_{2,j}(m),  \]  
\[ \widehat{a}^\ast_{1,j}(r) = < X_{1,j}^\ast-\overline{X}^\ast_{1,n_1}, \widehat{\varphi}^\ast_{r}>,  \quad  \widehat{a}^\ast_{2,j}(r) = < X_{2,j}^\ast-\overline{X}^\ast_{2,n_2}, \widehat{\varphi}^\ast_{r}> \]
and  $ \overline{X}^\ast_{i,n_i}(t)=n_i^{-1}\sum_{j=1}^{n_i} X_{i,j}^\ast(t)$, $ t \in {\cal I}$,  $ i =1,2$. Furthermore, 
  $ \widehat{\lambda}_{r}^\ast$ and $ \widehat{\varphi}^\ast_{r}$, $ r=1,2, \ldots, N$,  denote the eigenvalues and eigenfunctions, respectively, of the pooled covariance matrix 
 $ \widehat{C}_{N}^\ast(t,s) = (n_1/N)\widehat{C}^\ast_{1,n_1}(t,s)+  (n_2/N)\widehat{C}^\ast_{2,n_2}(t,s)$, where $ \widehat{C}^\ast_{i,n_1}=n_i^{-1}\sum_{j=1}^{n_i}(X_{i,j}^\ast(t) - \overline{X}^\ast_{i,n_i}(t))(X_{i,j}^\ast(s) - \overline{X}^\ast_{i,n_i}(s))$, $ t, s \in {\cal I}$, $i =1,2$.  We assume that 
  $ \widehat{\lambda}^\ast_1 \geq  \widehat{\lambda}^\ast_2 \geq  \widehat{\lambda}^\ast_3 \geq \cdots \geq \widehat{\lambda}^\ast_N$, and recall that $N=n_1+n_2$.  We first establish  the following  useful  lemmas.
  
\begin{lemma} \label{le.covlem1} Under the assumptions of Theorem~\ref{th:size-T2star} we have, conditionally on $ {\bf X}_N$, that,  for $ 1 \leq i \leq p$, 
$$
(i) \;\;|\widehat{\lambda}^\ast_i  - \widehat{\lambda}_i |  = O_P(N^{-1/2}) \quad \text{and}
\quad (ii)\;\;  \|\widehat{\varphi}^\ast_i -\widehat{c}^\ast_i \widehat{\varphi}_i\|  = O_P(N^{-1/2}),
$$
where $ \widehat{c}_i^\ast={\rm sign}(<\widehat{\varphi}^\ast_i,\widehat{\varphi}_i>)$.
\end{lemma}  
\noindent{\bf Proof:} \   We fist show that  for $ r \in \{1,2\}$, 
\begin{equation} \label{eq.lem1.1}
\Big\| n_r^{-1/2}\sum_{j=1}^{n_r}\big\{(X_{r,j}^\ast(t)-\overline{X}^\ast_{r,n_r}(t))(X_{r,j}^\ast(s)-\overline{X}^\ast_{r,n_r}(s)) - \widehat{C}_{N}(t,s)\big\} \Big\| = O_P(1)
\end{equation}
and
\begin{equation} \label{eq.lem1.2}
\Big\| n_r^{-1/2}\sum_{j=1}^{n_r}\big\{X_{r,j}^\ast(t)- \overline{X}^\ast_{r,n_r}(t)\big\} \Big\| = O_P(1).
\end{equation}
Let $ Y^\ast_{r,j}(t) = X_{r,j}^\ast(t)-\overline{X}^{\ast}_{r,n_r}(t)$ and  $ \widetilde{Y}^\ast_{r,j}(t) = X_{r,j}^\ast(t)-\overline{X}_{r,n_r}(t)$. Then
\begin{align*}
\EE\Big\|n_r^{-1/2}\sum_{j=1}^{n_r} \big\{Y^\ast_{r,j}(t)Y^\ast_{r,j}(s) -\widehat{C}_N(t,s)\big\}\Big\|^2 &
 \leq \frac{2}{n_r}\int_{0}^{1}\int_{0}^{1} \sum_{j=1}^{n_r} \EE\Big(\widetilde{Y}^\ast_{r,j}(t)\widetilde{Y}^\ast_{r,j}(s) - \widehat{C}_N(t,s) \Big  )^2  dt ds  \\
 & \ \ \ \ +
 2\; \EE\Big\|n_r^{-1/2}\sum_{j=1}^{n_r} \big\{Y^\ast_{r,j}(t)Y^\ast_{r,j}(s) -\widetilde{Y}^\ast_{r,j}(t)\widetilde{Y}^\ast_{r,j}(s) \big\}\Big\|^2 \\
 & = 2\int_{0}^{1}\int_{0}^{1} \EE\Big(\varepsilon^\ast_{r,1}(t)\varepsilon^\ast_{r,1}(s) - \widehat{C}_N(t,s) \Big)^2  dt ds + O_P(1),
\end{align*} 
using the fact that $Y^\ast_{r,j}(t) =\widetilde{Y}^\ast_{r,j}(t) +( \overline{X}_{r,n_r}(t)-\overline{X}^\ast_{r,n_r}(t)) = \varepsilon^\ast_{r,j}(t) + O_P(n_r^{-1/2})$. Now, since $\EE(\varepsilon^\ast_{r,1}(t)\varepsilon^\ast_{r,1}(s) - \widehat{C}_N(t,s))^2 = N^{-1}\sum_{r=1}^2\sum_{j=1}^{n_r}[\widehat{\varepsilon}_{r,j}(t)\widehat{\varepsilon}_{r,j}(s) - \widehat{C}_N(t,s)]^2=O_P(1)$, assertion (\ref{eq.lem1.1}) follows by Markov's inequality. Assertion (\ref{eq.lem1.2}) follows by the same inequality and because  
\begin{align*}
 \EE\Big\| n_r^{-1/2}\sum_{j=1}^{n_r}\big\{X_{r,j}^\ast(t)- \overline{X}^\ast_{r,n_r}(t)\big\} \Big\|^2   & \leq 2\; \EE\Big\| n_r^{-1/2}\sum_{j=1}^{n_r} \widetilde{Y}^\ast_{r,j} \Big\|^2 + O_P(1)\\
 &   = \int_{0}^{1}  \EE(\varepsilon^\ast_{r,1}(t))^2 dt + O_P(1) = O_P(1).\\
 \end{align*}  
 Using  (\ref{eq.lem1.1}), we get  that
 \begin{equation}
 \|\widehat{C}^\ast_N -\widehat{C}_N\|_{S}  \leq  \frac{n_1}{N}\|\widehat{C}^\ast_{1,n_1} - \widehat{C}_N\|_{S} +  \frac{n_2}{N}\|\widehat{C}^\ast_{2,n_2}-\widehat{C}_N\|_{S} = \frac{n_1}{N}O_P(n_1^{-1/2}) +  \frac{n_2}{N} O_P(n_2^{-1/2}) = O_P(N^{-1/2}).
 \label{eqFanis:123}
 \end{equation}
 Using (\ref{eqFanis:123}) and Lemmas 2.2 and 2.3 of  Horv\'ath \& Kokoszka (2012),  we have that, for $ 1 \leq i \leq p$, 
 \[ | \widehat{\lambda}^\ast_i -\widehat{\lambda}_i | \leq \| \widehat{C}_N^\ast - \widehat{C}_N\|_{S} =O_P(N^{-1/2}),\]
 and 
 \[ \| \widehat{\varphi}^\ast_i - \widehat{c}^\ast_i \widehat{\varphi}_i \| = O_P(\| \widehat{C}_N^\ast - \widehat{C}_N\|_{S}) = O_P(N^{-1/2}).\]
This completes the proof of the lemma.\hfill $\Box$

\bigskip

Let 
\[ \widetilde{A}^\ast_{1,n_1}(i,j) = \frac{1}{n_r}\sum_{k=1}^{n_r} <X_{1,k}^\ast-\overline{X}_{1,n_1}, \widehat{c}^\ast_i\widehat{\varphi}_i> <X_{1,k}^\ast-\overline{X}_{1,n_1}, \widehat{c}^\ast_j\widehat{\varphi}_j>,\]
\[ \widetilde{A}^\ast_{2,n_2}(i,j) = \frac{1}{n_2}\sum_{k=1}^{n_2} <X_{2,k}^\ast-\overline{X}_{2,n_2}, \widehat{c}^\ast_i\widehat{\varphi}_i> <X_{2,k}^\ast-\overline{X}_{2,n_2}, \widehat{c}^\ast_j\widehat{\varphi}_j>\]
and $ \widehat{\Delta}^\ast_N=( (\widehat{\Delta}^\ast_N(i,j)), 1 \leq i,j \leq p))$, where $ \widehat{\Delta}^\ast_N(i,j)=( \widetilde{A}^\ast_{1,n_1}(i,j) - \widetilde{A}^\ast_{2,n_2}(i,j)) $.

\begin{lemma} \label{le.covlem2} Under  the assumptions of Theorem~\ref{th:size-T2star} we have, conditionally on $ {\bf X}_N$, that 
\[ \sqrt{\frac{n_1n_2}{N}}\Big( \widehat{\Delta}^\ast_N  - \widetilde{\Delta}^\ast_N\Big) =  o_P(1).\]
\end{lemma}
 \noindent{\bf Proof:} \  Let $ \breve{a}^\ast_{r,j}(i) = < X_{r,j}^\ast - \overline{X}_{r,n_r}, \widehat{\varphi}^\ast_i>$, $ \widetilde{a}^\ast_{r,j}(i) = < X_{r,j}^\ast - \overline{X}_{r,n_r}, \widehat{c}_i^\ast \widehat{\varphi}_i>$. We first show that we can replace $ \widehat{a}^\ast_{r,j}(i)$ in $ \widehat{\Delta}^\ast_N(i,j)$  by 
 $ \breve{a}^\ast_{r,j}(i)$. For this, notice that,  for $ r \in \{1,2\}$, we have 
 \begin{align*}
 \sqrt{\frac{n_1n_2}{N}}\frac{1}{n_r}\sum_{k=1}^{n_r}\big(&  \widehat{a}^\ast_{r,k}(i)\widehat{a}^\ast_{r,k}(j) - \breve{a}^\ast_{r,k}(i)\breve{a}^\ast_{r,k}(j)\big) =  \sqrt{\frac{n_1n_2}{N}}\frac{1}{n_r}\sum_{k=1}^{n_r}\big(  \widehat{a}^\ast_{r,k}(i) - \breve{a}^\ast_{r,k}(i)\big)\widehat{a}^\ast_{r,k}(j)\\
&  -\sqrt{\frac{n_1n_2}{N}}\frac{1}{n_r}\sum_{k=1}^{n_r}\big(  \breve{a}^\ast_{r,k}(j) - \widehat{a}^\ast_{r,k}(j)\big)\breve{a}^\ast_{r,k}(i)
 \end{align*}
 and
 \begin{align*}
\sqrt{\frac{n_1n_2}{N}}\frac{1}{n_r}\sum_{k=1}^{n_r}\big(  \widehat{a}^\ast_{r,k}(i) - \breve{a}^\ast_{r,k}(i)\big)\widehat{a}^\ast_{r,k}(j) & = \sqrt{\frac{n_1n_2}{N}}\int_{0}^{1}\int_{0}^{1} (\overline{X}_{r,n_r}(t)-\overline{X}^\ast_{r,n_r}(t))\widehat{\varphi}^\ast_i(t)dt\\
& \times \frac{1}{n_r}\sum_{k=1}^{n_r}(X_{r,k}^\ast(s)-\overline{X}^\ast_{r,n_r}(s))\widehat{\varphi}^\ast_j(s)ds \\
&=0.
\end{align*}   
Let $ \breve{\Delta}^\ast_N(i,j)$ be the same expression as $ \widehat{\Delta}^\ast_N(i,j)$ with $\widehat{a}^\ast_{r,j}(l) $ replaced by $\breve{a}^\ast_{r,j}(l) $, $ l \in \{i,j\}$, and notice that by the previous considerations,  $\sqrt{n_1n_2/N}| \widehat{\Delta}^\ast_N(i,j)- \breve{\Delta}^\ast_N(i,j)|=o_P(1)$. Furthermore,  
\begin{align*}
\sqrt{\frac{n_1n_2}{N}}( \breve{\Delta}^\ast_N(i,j) -\widetilde{\Delta}^\ast_N(i,j)) & = \sqrt{\frac{n_1n_2}{N}}\int_{0}^{1}\int_{0}^{1} \frac{1}{n_1}\sum_{k=1}^{n_1}\Big[(X_{1,k}^\ast(t)-\overline{X}_{1,n_1}(t))(X_{1,k}^\ast(s)-\overline{X}_{1,n_1}(s))  - \widehat{C}_N(t,s)\Big]\\
& \times \Big(\widehat{\varphi}_i^\ast(t)\widehat{\varphi}^\ast_j(s) -\widehat{c}^\ast_i\widehat{c}_j^\ast  \widehat{\varphi}_i(t)\widehat{\varphi}_j(s) \Big)dtds\\
& -  \sqrt{\frac{n_1n_2}{N}}\int_{0}^{1}\int_{0}^{1} \frac{1}{n_2}\sum_{k=1}^{n_2}\Big[(X_{2,k}^\ast(t)-\overline{X}_{2,n_2}(t))(X_{2,k}^\ast(s)-\overline{X}_{2,n_2}(s))  - \widehat{C}_N(t,s)\Big]\\
& \times \Big(\widehat{\varphi}_i^\ast(t)\widehat{\varphi}^\ast_j(s) -\widehat{c}^\ast_i\widehat{c}_j^\ast  \widehat{\varphi}_i(t)\widehat{\varphi}_j(s) \Big)dtds\\
& = V_{1,N} + V_{2,N},
\end{align*}
with an obvious notation for $ V_{r,N}$, $ r \in \{1,2\}$.  Now, 
\begin{align*}
|V_{r,N}|  & \leq  \sqrt{\frac{n_1n_2}{N}}\int_{0}^{1}\int_{0}^{1} \Big|\frac{1}{n_r}\sum_{k=1}^{n_r}\Big[(X_{r,k}^\ast(t)-\overline{X}_{r,n_r}(t))(X_{r,k}^\ast(s)-\overline{X}_{r,n_r}(s))  - \widehat{C}_N(t,s)\Big]\\
& \ \ \ \  \times \Big(\widehat{\varphi}_i^\ast(t) -  \widehat{c}^\ast_i  \widehat{\varphi}_i(t)\Big)  \widehat{\varphi}^\ast_j(s)\Big|dt ds \\
& +    \sqrt{\frac{n_1n_2}{N}}\int_{0}^{1}\int_{0}^{1} \Big|\frac{1}{n_r}\sum_{k=1}^{n_r}\Big[(X_{r,k}^\ast(t)-\overline{X}_{r,n_r}(t))(X_{r,k}^\ast(s)-\overline{X}_{r,n_r}(s))  - \widehat{C}_N(t,s)\Big]\\
& \ \ \ \ \times \Big(\widehat{c}^\ast_j  \widehat{\varphi}_j(s) -   \widehat{\varphi}_j^\ast(t)\Big)  \widehat{c}^\ast_i \widehat{\varphi}_i(t)\Big|dt ds\\
& \leq \frac{1}{\sqrt{n_r}} \sqrt{\frac{n_1n_2}{N}}\Big\|\frac{1}{\sqrt{n_r}}\sum_{k=1}^{n_r} \Big[ (X_{r,k}^\ast(t)-\overline{X}_{r,n_r}(t))(X_{r,k}^\ast(s)-\overline{X}_{r,n_r}(s))  - \widehat{C}_N(t,s)\Big] \Big\|\\
& \ \ \ \ \times \Big\{\| \widehat{\varphi}^\ast_i - \widehat{c}_i^\ast \widehat{\varphi}_i\| + \| \widehat{\varphi}^\ast_j - \widehat{c}_j^\ast \widehat{\varphi}_j \| \Big\}\\
& = O_P\Big(\frac{1}{\sqrt{n_r}}  \sqrt{\frac{n_1n_2}{N}} \Big)\Big\{O_P(\| \widehat{\varphi}^\ast_i - \widehat{c}_i^\ast \widehat{\varphi}_i\| ) + O_P ( \| \widehat{\varphi}^\ast_j - \widehat{c}_j^\ast \widehat{\varphi}_j\|) \Big) = O_{P}( N^{-1/2}),
\end{align*}
because of (\ref{eq.lem1.1}) and Lemma~\ref{le.covlem1}.
This completes the proof of the lemma.  \hfill $\Box$
 
\bigskip
\noindent{\bf Proof of Theorem~\ref{th:size-TNstar}:} 
Let $ {\mathcal S} $ be the Hilbert space of Hilbert--Schmidt operators endowed with the inner product  $ \langle\Psi_1, \Psi_2 \rangle_{\mathcal S} =\sum_{j=1}^{\infty} \langle\Psi_1(e_j),\Psi_2(e_j) \rangle$ for $  \Psi_1$, $  \Psi_2 \in {\mathcal S}$, where  $ \{e_j:\ j=1,2,\ldots\}$ is an orthonormal basis in ${\mathcal H}$. Notice that  $ \widehat{\mathcal C}^\ast_i \in {\mathcal S} $, $ i=1,2$. Since  
$$ \overline{X}^\ast_{i,n_i} = \overline{X}_{i,n_i} + O_P(n_i^{-1/2}) \quad \text{and}\quad  n_i^{-1/2}\sum_{j=1}^{n_i}(X^\ast_{i,j}-\overline{X}^\ast_{i,n_i})=O_P(1), \quad i=1,2,$$  we get 
\begin{eqnarray*}
\widehat{\mathcal C}^\ast_i &=& \frac{1}{n_i}\sum_{j=1}^{n_i}(X^\ast_{i,j}-\overline{X}^\ast_{i,n_i}) \otimes (X^\ast_{i,j}-\overline{X}^\ast_{i,n_i}) \\ &=& 
\frac{1}{n_i}\sum_{j=1}^{n_i}(X^\ast_{i,j}-\overline{X}_{i,n_i}) \otimes (X^\ast_{i,j}-\overline{X}_{i,n_i}) + O_{P}(n^{-1}), \quad i=1,2,
\end{eqnarray*}
where  the random variables $(X^\ast_{i,j}-\overline{X}_{i,n_i}) \otimes (X^\ast_{i,j}-\overline{X}_{i,n_i})$ are, conditional on $X_N$, independent and identically distributed. By a central limit theorem for triangular arrays of independent and identically distributed ${\mathcal S}$-valued random variables (see, e.g., Politis \& Romano (1992, Theorem 4.2)), we get, conditionally on $X_N$,  that $ \sqrt{n_i}( \widehat{{\mathcal C}}^\ast_i - \widehat{{\mathcal C}}_N)$ converges weakly to a Gaussian random element ${\cal U}$ in $ {\mathcal S} $ with mean zero and covariance operator $ {\mathcal B}=\theta_1 {\mathcal B_1} +(1-\theta_1) {\mathcal B_2}$ as $ n_i\rightarrow \infty$. Here, 
$ {\mathcal B}_i$ is  the covariance operator of the limiting Gaussian random element $U_i$  to which $ \sqrt{n_i}(\widehat{\mathcal C}_i-{\mathcal C}_i)$ converges weakly as $ n_i\rightarrow \infty$.  

By the independence of the bootstrap random samples between the two populations, we have, conditional on $X_N$,
\begin{align*}
T_N^{\ast} & =N\|\widehat{\mathcal C}_1^\ast - \widehat{\mathcal C}_2^\ast\|^2_{\mathcal S} \\
& = N \langle
\widehat{\mathcal C}_1^\ast-  \widehat{\mathcal C}_N,
\widehat{\mathcal C}_1^\ast-  \widehat{\mathcal C}_N  \rangle_{\mathcal S}  + N \langle
\widehat{\mathcal C}_2^\ast-  \widehat{\mathcal C}_N,
\widehat{\mathcal C}_2^\ast-  \widehat{\mathcal C}_N  \rangle_{\mathcal S}  \\
&=  \frac{N}{n_1} \|\sqrt{n_1}(\widehat{\mathcal C}_1^\ast-\widehat{\mathcal C}_N)\|^2_{\mathcal{S}} + 
 \frac{N}{n_2}\|\sqrt{n_2}(\widehat{\mathcal C}_2^\ast-\widehat{\mathcal C}_N)\|^2_{\mathcal{S}}.
\end{align*}
Hence, taking into account the above results and that $ n_1/N \rightarrow \theta_1$, we have that $N\|\widehat{\mathcal C}_1^\ast - \widehat{\mathcal C}_2^\ast\|^2_{\mathcal S}$ converges weakly to $\sum_{l=1}^{\infty}\tilde{\lambda}_l Z_l^2$ as $ n_1, n_2\rightarrow \infty$, where $ \tilde{\lambda}_l$, $ l \geq 1$,  are the eigenvalues of the operator $ \widetilde{\mathcal B} = \theta_1^{-1}{\mathcal B}+ (1-\theta_1)^{-1}{\mathcal B}$ and $ Z_l$, $ l \geq 1$, are independent standard (real-valued) Gaussian distributed  random variables.  Since ${\mathcal B_1}={\mathcal B_2}$, the assertion follows. \hfill $\Box$

\bigskip
\noindent{\bf Proof of Theorem~\ref{th:size-T2star}:} \   Recall that  $ T_{p,N}^\ast=(n_1n_2/N)\widehat{\xi}^{\ast^{'}}_{N}\widehat{L}^{\ast^{-1}}_N \widehat{\xi}_{N}^\ast$ with 
$ \widehat{\xi}^\ast_N={\rm vech}(\widehat{\Delta}_N^\ast)$ and
  $\widehat{L}^{\ast}_N$ is an  estimator of the covariance matrix of $ \sqrt{n_1n_2/N}\widehat{\xi}^\ast_N$. The element of the latter matrix corresponding 
  to the covariance of $ \sqrt{n_1n_2/N}\widehat{\xi}^\ast_N(i_1,j_1) $ and $ \sqrt{n_1n_2/N} \widehat{\xi}^\ast_N(i_2,j_2)$, $ 1 \leq i_1\leq  j_1 \leq p$ and $  1 \leq i_2\leq  j_2 \leq p$,  is denoted by 
  $l(\widehat{\xi}^\ast_N(i_1,j_1), \widehat{\xi}^\ast_N(i_2,j_2)) $ and is estimated by
  \begin{align*}
 \widehat{l}(\widehat{\xi}^\ast_N(i_1,j_1), \widehat{\xi}^\ast_N(i_2,j_2)) & = \frac{n_2}{n_1+n_2}\Big\{ \frac{1}{n_1}\sum_{j=1}^{n_1}\widehat{a}^\ast_{1,j}(i_1)\widehat{a}^\ast_{1,j}(j_1)\widehat{a}^\ast_{1,j}(i_2)\widehat{a}^\ast_{1,j}(j_2) - <\widehat{C}^\ast_{1,n_1}\widehat{\varphi}^\ast_{i_1},\widehat{\varphi}^\ast_{j_1} >\\
 & \ \ \ \  \times   <\widehat{C}^\ast_{1,n_1}\widehat{\varphi}^\ast_{i_2},\widehat{\varphi}^\ast_{j_2} > \Big\} + 
  \frac{n_1}{n_1+n_2}\Big\{ \frac{1}{n_2}\sum_{j=1}^{n_2}\widehat{a}^\ast_{2,j}(i_1)\widehat{a}^\ast_{2,j}(j_1)\widehat{a}^\ast_{2,j}(i_2)\widehat{a}^\ast_{2,j}(j_2) \\
  & \ \ \ \  - <\widehat{C}^\ast_{2,n_2}\widehat{\varphi}^\ast_{i_1},\widehat{\varphi}^\ast_{j_1}> <\widehat{C}^\ast_{2,n_2}\widehat{\varphi}^\ast_{i_2},\widehat{\varphi}^\ast_{j_2} > \Big\}.
 \end{align*} 
Let  $\widehat{\Delta}_N(r,m)= \widehat{A}_{1,n_1}(r,m)-
\widehat{A}_{2,n_2}(r,m) $ be the $(r,m)$th element of the matrix $ \widehat{\Delta}_N^\ast$. To establish the theorem, it suffices to show that, under the assumptions made, 
the following assertions (\ref{eq.covth2part1}) and (\ref{eq.covth2part2}) are true.  
 \begin{equation} \label{eq.covth2part1}
 {\mathcal L}\Big( (\sqrt{\frac{n_1n_2}{N}}\widehat{\Delta}^\ast_N(i,j), 1 \leq i,j \leq p) \Big| {\bf X}_N\Big)  \Rightarrow {\mathcal L}( (\Delta(i,j), 1 \leq i, j \leq p)\Big),
 \end{equation}
where $ \Delta=(\Delta(i,j), 1 \leq i,j \leq p)$   is a Gaussian random matrix with $ \EE(\Delta(i,j))=0$, $ 1 \leq i,j,\leq p$ having a  positive definite covariance matrix $ \Sigma $ 
 with elements $ \sigma(i_1,j_1,i_2,j_2)=Cov( \Delta(i_1,j_1),\Delta(i_2,j_2) ) $,  $ 1\leq i_1,i_2,j_1,j_2 \leq p$, given by 
 \begin{align*}
 \sigma(i_1,j_1,i_2,j_2) =& (1-\theta) \Big\{ \EE(<X_{1,j}-\mu_1,\varphi_{i_1}><X_{1,j}-\mu_1,\varphi_{j_1}><X_{1,j}-\mu_1,\varphi_{i_2}><X_{1,j}-\mu_1,\varphi_{j_2}>)\\
 & \ \ - <C\varphi_{i_1},\varphi_{j_1}><C\varphi_{i_2},\varphi_{j_2}>\Big\}\\
   & \  \ +  \theta \Big\{ \EE(<X_{2,j}-\mu_2,\varphi_{i_1}><X_{2,j}-\mu_2,\varphi_{j_1}><X_{2,j}-\mu_2,\varphi_{i_2}><X_{2,j}-\mu_2,\varphi_{j_2}>)\\
 & \  \ - <C\varphi_{i_1},\varphi_{j_1}><C\varphi_{i_2},\varphi_{j_2}>\Big\},
 \end{align*}
 where $ C=(1-\theta)C_1 + \theta C_2$.
 Furthermore, 
 \begin{equation} \label{eq.covth2part2}
  \widehat{l}(\widehat{\xi}^\ast_N(i_1,j_1), \widehat{\xi}^\ast_N(i_2,j_2))  - 
  \widehat{c}^\ast_{i_1} \widehat{c}^\ast_{j_1} \widehat{c}^\ast_{i_2} \widehat{c}^\ast_{j_2}  \sigma(i_1,j_1,i_2,j_2) \rightarrow 0, \ \ \mbox{in probability}, 
 \end{equation}
 for all $ 1\leq i_1,i_2,j_1,j_2 \leq p$.  
 
 To establish (\ref{eq.covth2part1}), recall that by Lemma~\ref{le.covlem2} it suffices to consider the asymptotic distribution of $\sqrt{n_1n_2/N}\widetilde{\Delta}^\ast_N$. 
 Let $ \widetilde{Y}^\ast_{r,k}(i,j)=<X_{r,k}^\ast-\overline{X}_{r,n_r},\widehat{c}^\ast_i\widehat{\varphi}_i><X_{r,k}^\ast-\overline{X}_{r,n_r},\widehat{c}^\ast_j\widehat{\varphi}_j>$ and notice that 
 \[ \sqrt{n_1n_2/N}\widetilde{\Delta}^\ast_N (i,j)=  \sqrt{n_1n_2/N}\Big(\frac{1}{n_1}\sum_{k=1}^{n_1}\widetilde{Y}^\ast_{1,k}(i,j) -   \frac{1}{n_2}\sum_{k=1}^{n_2}\widetilde{Y}^\ast_{2,k}(i,j)\Big).\]
 Since $ \widehat{c}^\ast_{i} $  and $ \widehat{c}^\ast_{j}$ change solely the sign of $ \widetilde{Y}^\ast_{r,k}(i,j)$, they  do not affect the limiting distribution   of the two sums above.
 Thus, without loss of generality, we set $ \widehat{c}^\ast_i=\widehat{c}^\ast_j=1$. Let  
   $ \breve{Y}^\ast_{r,k}(i,j)=<X_{r,k}^\ast-\overline{X}_{r,n_r},\widehat{\varphi}_i><X_{r,k}^\ast-\overline{X}_{r,n_r},\widehat{\varphi}_j>$, notice that 
$ \EE(\breve{Y}^\ast_{r,k}(i,j))=  \widehat{\lambda}_{i} {\bf 1}_{\{i=j\}}$, and consider  instead of  the distribution of 
   $ \sqrt{n_1n_2/N}\widetilde{\Delta}^\ast_N$ the distribution of the asymptotically equivalent sum $  \sqrt{n_1n_2/N}Z^\ast_N(i,j)$, given by 
\begin{align*}
 \sqrt{n_1n_2/N} Z^\ast_N(i,j) & = \sqrt{n_1n_2/N}\Big(\frac{1}{n_1}\Big(\sum_{k=1}^{n_1}\breve{Y}^\ast_{1,k}(i,j) -  \frac{1}{n_2} \sum_{k=1}^{n_2}\breve{Y}^\ast_{2,k}(i,j)\Big) \\
 & = \sqrt{n_2/N}\frac{1}{\sqrt{n_1}}\sum_{k=1}^{n_1}(\breve{Y}^\ast_{1,k}(i,j) -  {\bf 1}_{\{i=j\}} \widehat{\lambda}_{i}) -   
 \sqrt{n_1/N}\frac{1}{\sqrt{n_2}}\sum_{k=1}^{n_2}(\breve{Y}^\ast_{2,k}(i,j)- {\bf 1}_{\{i=j\}} \widehat{\lambda}_{i})\\
 & = \sqrt{n_2/N} Z^\ast_{1,N}(i,j) - \sqrt{n_1/N} Z_{2,N}^\ast(i,j),
 \end{align*}
with an obvious notation for $Z^\ast_{r,N}(i,j)  $, $ r \in \{1,2\} $.
Notice that, conditionally on $ {\bf X}_N$, $Z^\ast_N(i,j)$ is distributed as the difference of the  two independent sums  $ Z^\ast_{1,N}(i,j)$ and $ Z^\ast_{2,N}(i,j)$,
where for $ r \in \{1,2\}$, $ Z^\ast_{r,N}(i,j)$ is a sum of  the independent and  identically distributed random variables $ \breve{Y}^\ast_{r,k}(i,j)-  \widehat{\lambda}_{i}{\bf 1}_{\{i=j\}}  $, $ k=1,2, \ldots, n_r$. Furthermore,   $\EE(Z^\ast_{r,N}(i,j) )  =0$ and since $ \varepsilon_{r,k}^\ast = X^\ast_{r,k} - \overline{X}_{r,n_r}$, we get   
\begin{align} \label{eq.sigma1}
 Cov(Z_{r,k}(i_1,j_1), Z_{r,k}(i_2,j_2))& = \EE (Z_{r,k}(i_1,j_1) Z_{r,k}(i_2,j_2) ) \nonumber \\
 & = \EE\Big[<\varepsilon^\ast_{1,k},\widehat{\varphi}_{i_1}><\varepsilon^\ast_{1,k},\widehat{\varphi}_{j_1}><\varepsilon^\ast_{1,k},\widehat{\varphi}_{i_2}><\varepsilon^\ast_{1,k},\widehat{\varphi}_{j_2}> \Big] \nonumber \\
 & \ \ \  \  -\widehat{\lambda}_{i_1} \widehat{\lambda}_{i_2}{\bf 1}_{\{  i_1 =j_1\}}{\bf 1}_{\{  i_2 =j_2\}} \nonumber \\
 & = \sum_{r=1}^{2}\frac{n_r}{N}\frac{1}{n_r}\sum_{k=1}^{n_r}<\widehat{\varepsilon}_{r,k},\widehat{\varphi}_{i_1}><\widehat{\varepsilon}_{r,k},\widehat{\varphi}_{j_1}><\widehat{\varepsilon}_{r,k},\widehat{\varphi}_{i_2}><\widehat{\varepsilon}_{r,k},\widehat{\varphi}_{j_2}> \nonumber \\
 & \ \ \ \  \  -\widehat{\lambda}_{i_1} \widehat{\lambda}_{i_2}{\bf 1}_{\{  i_1 =j_1\}}{\bf 1}_{\{  i_2 =j_2\}} \nonumber \\
 & \rightarrow (1-\theta)\; \EE\Big[<X_{1,k}-\mu_1,\varphi_{i_1}> <X_{1,k}-\mu_1,\varphi_{j_1}> <X_{1,k}-\mu_1,\varphi_{i_2}> \nonumber \\
 & \ \ \ \ \times  <X_{1,k}-\mu_1,\varphi_{j_2}>\Big]  + \theta\;  \EE\Big[<X_{2,k}-\mu_2,\varphi_{i_1}> <X_{2,k}-\mu_2,\varphi_{j_1}> \nonumber \\
 & \ \ \ \ \times <X_{2,k}-\mu_2,\varphi_{i_2}> <X_{2,k}-\mu_2,\varphi_{j_2}>\Big] -\lambda_{i_1} \lambda_{i_2}{\bf 1}_{\{  i_1 =j_1\}}{\bf 1}_{\{  i_2 =j_2\}}\nonumber \\
 & = \sigma(i_1,j_1,i_2,j_2),  
\end{align}
in probability, by the weak law of large numbers and using  $ <\widehat{\varepsilon}_{r,l},\widehat{\varphi}_i>= <X_{r,l}-\overline{X}_{r,n_r},\widehat{\varphi}_i>$.
Thus, by a multivariate central limit theorem for triangular arrays of real valued random vectors, we get that 
$ {\mathcal L}(  Z^\ast_{r,N} ) \Rightarrow Z$,  where $ Z$ is a Gaussian distributed $ p\times p$  random matrix, with $ \EE(Z(i,j))=0$ and $ Cov(Z(i_1,j_1),Z(i_2,j_2))=\sigma(i_1,j_1,i_2,j_2)$. 
To conclude the proof of (\ref{eq.covth2part1}), notice that 
\[ {\mathcal L}( Z^\ast_N) = {\mathcal L}(\sqrt{n_2/N} Z^\ast_{1,N} + \sqrt{n_1/N}Z^\ast_{2,N}) \Rightarrow {\mathcal L}(\sqrt{1-\theta} Z_1 + \sqrt{\theta} Z_2) = {\mathcal L}(Z), \]
where  $Z_1$ and $ Z_2$ are two independent copies of  the Gaussian random matrix $ Z$.   

To establish  (\ref{eq.covth2part2}), notice that, for $ r \in \{1,2\}$,
\begin{align*} 
\frac{1}{n_r}\sum_{j=1}^{n_r}& \widehat{a}^\ast_{r,j}(i_1) \widehat{a}^\ast_{r,j}(j_1)\widehat{a}^\ast_{r,j}(i_2)\widehat{a}^\ast_{r,j}(j_2) = \frac{1}{n_r}\sum_{j=1}^{n_r}\widetilde{a}^\ast_{r,j}(i_1)\widetilde{a}^\ast_{r,j}(j_1)\widetilde{a}^\ast_{r,j}(i_2)\widetilde{a}^\ast_{r,j}(j_2) + O_P(n_r^{-1/2}) \nonumber 
\end{align*}
and that,  for 
\begin{align*}
   \sigma^{(1)}(i_1,j_1,i_2,j_2) & = (1-\theta)\; \EE[ <X_{1,k}-\mu_1,\varphi_{i_1}> <X_{1,k}-\mu_1,\varphi_{j_1}> <X_{1,k}-\mu_1,\varphi_{i_2}>  <X_{1,k}-\mu_1,\varphi_{j_2}>] \\
  &   + \theta\; \EE[ <X_{2,k}-\mu_2,\varphi_{i_1}> <X_{2,k}-\mu_2,\varphi_{j_1}> <X_{2,k}-\mu_2,\varphi_{i_2}> <X_{2,k}-\mu_2,\varphi_{j_2}>],
\end{align*}
we have 
\begin{align}  \label{eq.l1}
 \frac{1}{n_r}\sum_{j=1}^{n_r}\widetilde{a}^\ast_{r,j}(i_1)\widetilde{a}^\ast_{r,j}(j_1)\widetilde{a}^\ast_{r,j}(i_2)\widetilde{a}^\ast_{r,j}(j_2) -   \widehat{c}^\ast_{i_1} \widehat{c}^\ast_{j_1} \widehat{c}^\ast_{i_2} \widehat{c}^\ast_{j_2}   \sigma^{(1)}(i_1,j_1,i_2,j_2)  \rightarrow 0,
\end{align}
in probability, since  as in  obtaining (\ref{eq.sigma1}), 
\begin{align*}
\EE\Big(\frac{1}{n_r}\sum_{j=1}^{n_r}& \widetilde{a}^\ast_{r,j}(i_1)\widetilde{a}^\ast_{r,j}(j_1)\widetilde{a}^\ast_{r,j}(i_2)\widetilde{a}^\ast_{r,j}(j_2)\Big) - 
\widehat{c}^\ast_{i_1} \widehat{c}^\ast_{j_1} \widehat{c}^\ast_{i_2} \widehat{c}^\ast_{j_2}   \sigma^{(1)}(i_1,j_1,i_2,j_2)  \\
& = \EE\Big[<\varepsilon^\ast_{r,k},\widehat{c}^\ast_{i_1}\widehat{\varphi}_{i_1}><\varepsilon^\ast_{r,k},\widehat{c}^\ast_{j_1}\widehat{\varphi}_{j_1}><\varepsilon^\ast_{r,k},\widehat{c}^\ast_{i_2}\widehat{\varphi}_{i_2}><\varepsilon^\ast_{r,k},\widehat{c}^\ast_{j_2}\widehat{\varphi}_{j_2}> \Big]  \\
&  \ \ \ \ -  \widehat{c}^\ast_{i_1} \widehat{c}^\ast_{j_1} \widehat{c}^\ast_{i_2} \widehat{c}^\ast_{j_2}   \sigma^{(1)}(i_1,j_1,i_2,j_2)  \\
& \rightarrow 0, 
\end{align*}
in probability, and also $ Var( n_r^{-1}\sum_{j=1}^{n_r}\widetilde{a}^\ast_{r,j}(i_1)\widetilde{a}^\ast_{r,j}(j_1)\widetilde{a}^\ast_{r,j}(i_2)\widetilde{a}^\ast_{r,j}(j_2) ) =O_P(n_r^{-1}),$ due to the independence of the random variables $\widetilde{a}^\ast_{r,j}(i_1)\widetilde{a}^\ast_{r,j}(j_1)\widetilde{a}^\ast_{r,j}(i_2)\widetilde{a}^\ast_{r,j}(j_2)$, for different $j$'s.
Furthermore,  by the triangular inequality and because  
$| <\widehat{C}_N\widehat{\varphi}_i,\widehat{\varphi}_j >- <C\varphi_i,\varphi_j > |\rightarrow 0$ in probability, it yields that 
\begin{equation} \label{eq.l2}
 |<\widehat{C}^\ast_{r,n_r}\widehat{\varphi}^\ast_i, \widehat{\varphi}^\ast_j> - \widehat{c}^\ast_i \widehat{c}^\ast_j <C\varphi_i,\varphi_j >| \rightarrow 0,
\end{equation}
in probability, 
 since  
\begin{align*}
|<\widehat{C}^\ast_{r,n_r}\widehat{\varphi}^\ast_i, \widehat{\varphi}^\ast_j> - \widehat{c}^\ast_i\widehat{c}^\ast_j<\widehat{C}_N\widehat{\varphi}_i,\widehat{\varphi}_j >| & \leq
 |<(\widehat{C}^\ast_{r,n_r}-\widehat{C}_N)\widehat{\varphi}^\ast_i, \widehat{\varphi}^\ast_j> |+ |<(\widehat{C}_{N}(\widehat{\varphi}^\ast_i-\widehat{c}^\ast_i\widehat{\varphi}_i), \widehat{\varphi}^\ast_j> | \\
 & \ \ \ \ +  |<(\widehat{C}_{N}\widehat{\varphi}_i, \widehat{\varphi}^\ast_j-\widehat{c}^\ast_j\widehat{\varphi}_j> |\\
 & = O_P\Big(\| \widehat{C}^\ast_{r,n_r}-\widehat{C}_N\| +
  \|  \widehat{\varphi}^\ast_i-\widehat{c}^\ast_i\widehat{\varphi}_i\| + \| \widehat{\varphi}^\ast_j-\widehat{c}^\ast_j\widehat{\varphi}_j \|\Big) \rightarrow 0,
\end{align*}
by  (\ref{eq.lem1.1}) and Lemma~\ref{le.covlem1}.   Equations (\ref{eq.l1}) and (\ref{eq.l2}) imply then assertion (\ref{eq.covth2part2}). 
This completes the proof of the theorem.  \hfill $\Box$
 
 \bigskip
 
\noindent{\bf Proof of Theorem~\ref{th:size-T1star}:} \   
Notice that under Gaussianity  of  the random functions $X_{1,j}$ and $ X_{2,j}$, the random variables $ <X_{1,j}-\mu_1,\varphi_i>$ and   $<X_{2,j}-\mu_2,\varphi_i>$ are independent Gaussian distributed with mean zero and variance $ \lambda_i$, $i=1,2,\ldots,p$. From assertion  (\ref{eq.covth2part1})  in the proof of Theorem~\ref{th:size-T2star}, we get that 
in this case,  the random variables $ \Delta(i,j)$ are for $ 1 \leq i\leq j \leq p$ independent with mean zero and $ Var(\Delta(i,j))=2\lambda_i^2$ if $i=j$ and $ Var(\Delta(i,j))=\lambda_i\lambda_j$ if $ i\neq j$. We then have that 
\begin{align*}
T_{p,N}^{\ast (G)} & = \frac{n_1n_2}{N}\frac{1}{2}\Big(\sum_{r=1}^p \frac{\widehat{\Delta}^{\ast^2}(r,r)}{\widehat{\lambda}^{\ast^2}_r} + 2\sum_{1\leq r < m \leq p} \frac{\widehat{\Delta}^{\ast^2}(r,m)}{\widehat{\lambda}^{\ast}_r \widehat{\lambda}^{\ast}_m} \Big)\\
& = \sum_{r=1}^p\Big(\sqrt{ \frac{n_1n_2}{N}}
 \frac{\widehat{\Delta}^{\ast}(r,r)}{\sqrt{2}\lambda_r}\Big)^2\Big( \frac{\lambda_r}{\widehat{\lambda}_m}\Big)^2 + 
 \sum_{1\leq r < m \leq p}\Big(\sqrt{\frac{n_1n_2}{N}}\frac{\widehat{\Delta}^{\ast}(r,m)}{\lambda_r \lambda_m} \Big)^2\Big( \frac{\lambda_r\lambda_m}{
 \widehat{\lambda}_r\widehat{\lambda}_m}\Big)^2\\
 & \Rightarrow \chi^2_{p(p+1)/2},
\end{align*}
since  by assertion (\ref{eq.covth2part1})  we have that 
$ \sqrt{n_1n_2/N}\widehat{\Delta}^{\ast}(r,r)/(\sqrt{2}\lambda_r)$ resp. $\sqrt{n_1n_2/N}\widehat{\Delta}^{\ast}(r,m)/(\lambda_r \lambda_m)  $ are asymptotically independent, standard Gaussian distributed random variables, and,  by Lemma 3   of  Fremdt {\em et al.} (2012)  and Lemma~\ref{le.covlem1}(i),  we get that 
$ \widehat{\lambda}_i^\ast \rightarrow \lambda_i$, in probability, for $ i=1,2, \ldots, p$. 
This completes the proof of the theorem. \hfill $\Box $
  
  \bigskip
  
\noindent{\bf Proof of Theorem~\ref{th:size-Ustar}:} \   
Define 
$$
Z^+_{n_1,n_2}(t) =\Big[ n_1^{-1/2}\sum_{j=1}^{n_1}\{X_{1,j}^+(t) - \overline{X}_{N}(t)\}, n_2^{-1/2}\sum_{j=1}^{n_2}\{X_{2,j}^+(t)
 - \overline{X}_{N}(t)\} \Big], \ t \in {\cal I},
 $$
 and
 $$
 Z^+_{i,n_i}(t)=n_i^{-1/2}\sum_{j=1}^{n_i}\{X_{i,j}^+(t) - \overline{X}_{N}(t)\}, \ t \in {\cal I}, \ \ i=1,2.
 $$
 Notice that, conditionally on $X_N$, $ Z^+_{1,n_1}(t) $ and $ Z^+_{2,n_2}(t) $ are independent, have covariance operators $ \widehat{\mathcal{C}}_{1}$ and $ \widehat{\mathcal{C}}_{2}$, respectively, and
 $X_{1,j}^+(t) $ and $X_{2,j}^+(t) $ have the same mean function $ \overline{X}_N(t)$.
 By a central limit theorem  for triangular arrays of independent and identically distributed
 $\mathcal{H}$-valued random variables (see, e.g., Politis \& Romano (1992, Theorem 4.2)),  it follows that, conditionally on $X_N$, $Z_{i,n_i}^+$ converges weakly to  a Gaussian random element ${\cal U}_i$ with mean zero and covariance operator $\mathcal{C}_i$ as $ n_i \rightarrow \infty$. 

By the independence of $  Z^+_{1,n_1} $ and $ Z^+_{2,n_2} $, we have, conditionally on $X_N$,  
\begin{align*}
S_N^+ & = \frac{n_1n_2}{N}\int_0^1 \{\overline{X}_{1,n_1}^+(t) - \overline{X}_{2,n_2}^+(t)\}^2dt\\
 & =   \frac{n_1n_2}{N}\int_0^1 \Big[\frac{1}{n_1}\sum_{t=1}^{n_1}\{X_{1,j}^+(t)-\overline{X}_{N}(t)\} - \frac{1}{n_2}\sum_{t=1}^{n_2}\{X_{2,j}^+(t)-\overline{X}_{N}(t)\} \Big]^2dt \\
 & = \int_0^1 \Bigg[\sqrt{\frac{n_2}{N}} Z_{1,n_1}^+(t)  - \sqrt{\frac{n_1}{N}} Z_{2,n_2}^+(t) \Bigg]^2dt,
 \end{align*}
 from which, and taking into account that $ n_1/N \rightarrow \theta_1$, we have that $S_N^+$ converges weakly to $\int_0^1 \Gamma^2(t) dt$ as $n_1, n_2\rightarrow \infty$. Hence, the assertion follows. \hfill $\Box$
 
 \bigskip
 
 To prove Theorem~\ref{th:size-S1star},  we first fix some notation. Let 
 \[ \widehat{a}^+_i=<\overline{X}^+_{1,n_1} - \overline{X}^+_{2,n_2},\widehat{\phi}_i^+> = \int_0^1(\overline{X}^+_{1,n_1}(t) - \overline{X}^+_{2,n_2}(t))\widehat{\phi}^+_i(t)dt,\ \ i=1,2, \ldots, p,\]
where $ \widehat{\phi}_i$, $ i=1,2, \ldots, p$ are the eigefunctions of $ \widehat{C}^+_N =n_2/N \widehat{C}^+_{1,n_1} + n_1/N \widehat{C}^+_{2,n_2}$ with 
$ \widehat{C}^+_{1,n_1}(t,s) =n_1^{-1}\sum_{j=1}^{n_1}(X^+_{1,j}(t)-\overline{X}_{1,n_1}(t)) (X^+_{1,j}(s)-\overline{X}_{1,n_1}(s))  $ and $ \widehat{C}^+_{2,n_2}(t,s) =n_2^{-1}\sum_{j=1}^{n_2}(X^+_{2,j}(t)-\overline{X}_{2,n_2}(t)) (X^+_{2,j}(s)-\overline{X}_{2,n_2}(s))   $.   Let $ \widehat{\tau}_i^+$ 
 be the eigenvalues corresponding to the eigenfunctions $ \widehat{\phi}_i^+$, $ i=1,2, \ldots, p$,  of $   \widehat{C}^+_N$. \\
 
The  following lemma is proved along the same lines as Lemma~\ref{le.covlem1} and is useful in establishing  Theorem~\ref{th:size-S1star}. Hence, its proof is omitted.

\begin{lemma} \label{le.meanlem1} Under the assumptions of Theorem~\ref{th:size-S1star} we have, conditionally on $ {\bf X}_N$ that,  for $ 1 \leq i \leq p$, 
$$
(i)\;\; |\widehat{\tau}^+_i  - \widehat{\tau}_i |  = O_P(N^{-1/2}) \quad \text{and} \quad
(ii)\;\; \|\widehat{\phi}^+_i -\widehat{c}^+_i \widehat{\phi}_i\|  = O_P(N^{-1/2}),
$$
where $ \widehat{c}_i^+={\rm sign}(<\widehat{\phi}^+_i,\widehat{\phi}_i>)$.
\end{lemma}    

\bigskip

 \noindent{\bf Proof of Theorem~\ref{th:size-S1star}:} \  
 Let $ \widehat{a}^+(1,p)=(\widehat{a}^+_1, \widehat{a}^+_2, \ldots,  \widehat{a}^+_p)^{'}$. We first show that 
 \begin{equation} \label{eq.asynormean}
 {\mathcal L} \Big(\sqrt{\frac{n_1n_2}{N}} \widehat{a}^+(1,p) \Big| {\bf X}_N\Big) \Rightarrow N(0, T), 
 \end{equation}
 where $ T=(T(i,j))_{i,j=1,2, \ldots, p}$ is a $p\times p$ diagonal matrix with $ T(i,i)=\tau_i$.   Let 
 $\widetilde{a}^+_i=<\overline{X}^+_{1,n_1} - \overline{X}^+_{2,n_2},\widehat{c}^+_i\widehat{\phi}_i> $ and $ 
   \widetilde{a}^+(1,p)=(\widetilde{a}^+_1, \widetilde{a}^+_2, \ldots,  \widetilde{a}^+_p)^{'} $.  We have that 
\begin{align*}
\sqrt{\frac{n_1n_2}{N}}|\widehat{a}^+_i -\widetilde{a}^+_i|  & =  \sqrt{\frac{n_1n_2}{N}}\big|<\overline{X}_{1,n_1} - \overline{X}_{2,n_2},\widehat{\phi}_i^+-\widehat{c}^+_i\widehat{\phi}_i>|\\
 & \leq \|\widehat{\phi}_i^+-\widehat{c}_i^+\widehat{\phi}_i\|  \Big\|  \sqrt{\frac{n_1n_2}{N}}\big(\overline{X}_{1,n_1} - \overline{X}_{2,n_2}\big)\Big\|\\
 &=O_P( \|\widehat{\phi}_i^+-\widehat{c}_i^+\widehat{\phi}_i\|  )\rightarrow 0,
\end{align*}  
 by Lemma~\ref{le.meanlem1}(ii).  Thus $ \sqrt{n_1n_2/N}\,\widehat{a}^+(1,p)=\sqrt{n_1n_2/N}\, \widetilde{a}^+(1,p)+o_P(1)$.  Now, let  for $ r \in \{1,2\}$, 
 $ L_{r,j}^+(k)=<X_{r,j}^+-\overline{X}_N,\widehat{c}_k^+\widehat{\phi}_k>$ and notice that   
 \begin{align*}
 \sqrt{n_1n_2/N}\, \widetilde{a}^+(1,p)& = \Big((\sqrt{\frac{n_2}{N}} \frac{1}{\sqrt{n_1}}\sum_{j=1}^{n_1} L^+_{1,j}(k) -
  \sqrt{\frac{n_1}{N}}\frac{1}{\sqrt{n_2}}\sum_{j=1}^{n_2} L^+_{2,j}(k) \big), k=1,2, \ldots, p \Big).
 \end{align*}
 Furthermore,  the $ L_{r,j}^+(k)$'s are independent  and satisfy $ E( L^+_{r,j}(k))=0$ and 
 \begin{align*}
 Cov\Big(n_r^{-1/2}\sum_{j=1}^{n_r} L^+_{r,j}(k_1), n_r^{-1/2}\sum_{j=1}^{n_r} L^+_{r,j}(k_2)\Big)& = \EE\Big(L^+_{r,1}(k_1) L^+_{r,1}(k_2)) \Big)\\
 &= \widehat{c}^+_{k_1}\widehat{c}^+_{k_2}\int_0^1\int_0^1 \EE(\varepsilon_{r,j}^+(t)\varepsilon_{r,j}(s))\widehat{\phi}_{k_1}(t)\widehat{\phi}_{k_2}(s)dtds\\
 &=  \widehat{c}^+_{k_1}\widehat{c}^+_{k_2}\int_0^1\int_0^1 \widehat{C}_{n,r}(t,s) \widehat{\phi}_{k_1}(t)\widehat{\phi}_{k_2}(s)dtds.
 \end{align*}
 This implies that  $ E( \sqrt{n_1n_2/N} \widetilde{a}^+(1,p))=0$ and that 
 \begin{align*}
 Cov \Big( \sqrt{n_1n_2/N} \widetilde{a}^+_{k_1}, \sqrt{n_1n_2/N} \widetilde{a}^+_{k_2} \Big) &  = 
 \widehat{c}^+_{k_1}\widehat{c}^+_{k_2} \int_0^1\int_0^1\Big( n_2/N\widehat{C}_{1,n_1}(t,s) + n_1/N \widehat{C}_{2,n_2}(t,s)\Big)\widehat{\phi}_{k_1}(t)
 \widehat{\phi}_{k_2}(s)dtds\\
 & = \widehat{c}^+_{k_1}\widehat{c}^+_{k_2} \int_0^1\int_0^1\widehat{C}_N(t,s) \widehat{\phi}_{k_1}(t)
 \widehat{\phi}_{k_2}(s)dtds\\
 & = {\bf 1}_{\{k_1=k_2\}}\widehat{\tau}_{k_1} \rightarrow \tau_{k_1}. 
 \end{align*}
Hence, (\ref{eq.asynormean}) follows then  by a multivariate central limit theorem  for triangular arrays of independent random variables.

Now,  (17) and   Lemma~\ref{le.meanlem1}(i) lead to  assertion (i) of the theorem,  since  
\begin{align*}
S_{p,N}^{+ (1)} & =\frac{n_1n_2}{N} \sum_{k=1}^p (\widehat{a}^+_k)^2/\widehat{\tau}_{k} - \sum_{k=1}^{p}\frac{\widehat{\tau}_k^+ \ - \ \widehat{\tau}_k}{\widehat{\tau}^+_k}\Big(\sqrt{n_1n_2/N}\widehat{a}_k^+/\widehat{\tau}_k\Big)^2\\
& =  \frac{n_1n_2}{N} \sum_{k=1}^p (\widehat{a}^+_k)^2/\widehat{\tau}_{k} + O_P\Big(\max_{1\leq k\leq p}|\widehat{\tau}^+_k-\widehat{\tau}_k| \Big)  \Rightarrow \chi^2_p, 
\end{align*}
and to  assertion (ii), since
\begin{align*}
S_{p,N}^{+ (2)} & = \sum_{k=1}^{p}\widehat{\tau}_k\Big(\sqrt{n_1n_2/N}\frac{\widehat{a}_{k}^+}{\sqrt{\widehat{\tau}k}} \Big)^2 \Rightarrow \sum_{k=1}^p\tau_k N^2_k,
\end{align*}
where $ N_k$, $k=1,2, \ldots, p$,  are independent, standard Gaussian random variables. 
This completes the proof of the theorem.  \hfill $\Box$

\section*{Acknowledgements}
We would like to thank Dr. Stefan Fremdt for providing us with the medfly data and the R-codes used in Fremdt {\em et al.} (2012).

\newpage

\renewcommand{\baselinestretch} {0.86}
\small\normalsize

\begin{table}[]
\begin{tabular}{|cccc|ccccccc|}
\hline
\hline
                   &              &         &  &  & NG& &&&BM &  \\
       $n_1$ &$ n_2$ & $p$ & Test-Stat. & $\alpha=0.01$ & $\alpha=0.05$ & $\alpha=0.10$ & & $\alpha=0.01$ & $\alpha=0.05$ & $\alpha=0.10$\\
\hline
25 & 25 & 2  &  $T_{p,N}$-Asym & 0.002 & 0.010& 0.040 && 0.000 & 0.024 & 0.058\\
         &        &    & $T_{p,N}$-Boot & 0.008& 0.048 & 0.104 && 0.004 & 0.054 & 0.116 \\ 
     &        &    &   $T_{p,N}^{(G)}$-Asym & 0.072 & 0.188& 0.292 && 0.008 & 0.048 & 0.100\\
              &        &    &   $T_{p,N}^{(G)}$-Boot & 0.006 & 0.040& 0.102 && 0.010 & 0.054  & 0.117\\
         25 & 25 & 3  &  $T_{p,N}$-Asym & 0.000 & 0.010& 0.034 && 0.000 & 0.020 & 0.048\\
         &        &    &       $T_{p,N}$-Boot & 0.006& 0.046 & 0.096 && 0.003 & 0.032 & 0.068\\ 
     &        &    & $T_{p,N}^{(G)}$-Asym & 0.068 & 0.198& 0.320 && 0.006 & 0.048 & 0.106\\
                        &        &    & $T_{p,N}^{(G)}$-Boot & 0.008 & 0.046& 0.094 && 0.004 & 0.038 & 0.088\\
              25 &  25 &-- & $T_{N}$-Boot &0.006 &0.034 &0.112  && 0.004  & 0.060  & 0.134\\
              50 & 50 & 2  &  $T_{p,N}$-Asym & 0.001 & 0.020& 0.054 && 0.002 & 0.026 & 0.064\\
              &        &    &    $T_{p,N}$-Boot & 0.010& 0.042 & 0.082 && 0.006 & 0.056 & 0.108\\
      &       &    &  $T_{p,N}^{(G)}$-Asym & 0.086 & 0.244& 0.328 && 0.008 & 0.050 & 0.112\\
                 &        &    & $T_{p,N}^{(G)}$-Boot & 0.004 & 0.069& 0.118 && 0.006 & 0.052 & 0.108 \\ 
         50 & 50 & 3  & $T_{p,N}$-Asym & 0.002 & 0.012& 0.050  && 0.000 & 0.016 & 0.046\\
         &        &    &   $T_{p,N}$-Boot & 0.006& 0.052 & 0.092 && 0.012 & 0.042 & 0.094\\ 
     &        &    &   $T_{p,N}^{(G)}$-Asym & 0.124 & 0.254& 0.340 && 0.006 & 0.048 & 0.094\\
                        &        &    &    $T_{p,N}^{(G)}$-Boot & 0.004 & 0.040& 0.114 &&  0.008& 0.042& 0.093\\
                        50 &  50 &-- & $T_{N}$-Boot &0.010 &0.056 &0.088  &&0.016   & 0.052 &0.110 \\
              100 & 100 & 2  &  $T_{p,N}$-Asym & 0.000 & 0.018& 0.040 && 0.002 & 0.022 & 0.048\\
              &        &    &    $T_{p,N}$-Boot & 0.004& 0.046 & 0.100 && 0.002 & 0.030 & 0.082\\
      &       &    &  $T_{p,N}^{(G)}$-Asym & 0.128 & 0.272& 0.376 && 0.006 & 0.042 & 0.108\\             
          &        &    &  $T_{p,N}^{(G)}$-Boot & 0.006 & 0.042& 0.090 && 0.002 & 0.022 & 0.074\\
         100 & 100 & 3  & $T_{p,N}$-Asym & 0.004 & 0.018& 0.044  && 0.002 & 0.034 & 0.060\\ 
         &        &    &   $T_{p,N}$-Boot & 0.006& 0.054& 0.094 &&0.002 & 0.037& 0.075\\
     &        &    &   $T_{p,N}^{(G)}$-Asym & 0.146 & 0.312& 0.410 && 0.008 & 0.048 & 0.090\\                 
          &        &    &   $T_{p,N}^{(G)}$-Boot & 0.006& 0.042& 0.100 && 0.009& 0.050&0.092\\
          100 &  100 &-- & $T_{N}$-Boot &0.006 &0.030 &0.096  && 0.006   & 0.052 &0.108 \\
\hline\hline
\end{tabular}
\caption{Empirical size of the tests for the equality of two covariance functions, based on the statistics $T_{p,N}^{(G)}$-Asym and $T_{p,N}$-Asym (asymptotic approximations) and $T_{p,N}^{(G)}$-Boot, $T_{p,N}$-Boot and $T_N$-Boot (bootstrap approximations), using two ($p=2$) and three ($p=3$) FPC's, both for Gaussian (BM) and non-Gaussian (NG) data. The curves in each sample were generated according to Brownian motions for Gaussian data and according to (\ref{eq:NG-sim}) for non-Gaussian data.}
\label{tab:empsize-cov}
\end{table}

\begin{table}[]
\centering
\begin{tabular}{|cc|ccccccc|}
\hline
\hline
                          &  &  & N=M=25& &&&N=M=50 &  \\
       $\gamma$ & Test-Stat. & $\alpha=0.01$ & $\alpha=0.05$ & $\alpha=0.10$ & & $\alpha=0.01$ & $\alpha=0.05$ & $\alpha=0.10$\\
\hline
2.0 & $T_{2,N}$-Asym & 0.000 & 0.042& 0.220 & &0.074  & 0.610 &0.850 \\
        & $T_{2,N}$-Boot & 0.008& 0.188 & 0.448 && 0.398 & 0.834 &  0.940\\     
        & $T_{N}$-Boot & 0.256& 0.622 & 0.784 && 0.574 & 0.860 &  0.930\\
         &  &  &  & &  & & &\\
2.2    &   $T_{2,N}$-Asym & 0.000 & 0.058& 0.318 &&0.196  & 0.774 &0.932 \\
              &   $T_{2,N}$-Boot & 0.022 & 0.328& 0.594 && 0.696 & 0.938  & 0.980\\
              & $T_{N}$-Boot & 0.332& 0.718 & 0.864 && 0.702 & 0.904 &  0.950\\
              &  &  &  & &  & & &\\
2.4    &   $T_{2,N}$-Asym & 0.000 & 0.132& 0.440 &&  0.358& 0.902 & 0.976\\
              &   $T_{2,N}$-Boot & 0.044 & 0.460& 0.732 &&  0.818&0.976   &0.990 \\
              & $T_{N}$-Boot & 0.382& 0.800 & 0.918 && 0.714 & 0.932 &  0.980\\  
              &  &  &  & &  & & &\\  
2.6    &   $T_{2,N}$-Asym & 0.000 & 0.162& 0.532 && 0.500 & 0.946 & 0.988\\
              &   $T_{2,N}$-Boot & 0.082 & 0.596& 0.840 && 0.900 & 0.990  & 0.996\\
              & $T_{N}$-Boot & 0.460& 0.810 & 0.914 &&  0.806&  0.956&  0.982\\   
              &  &  &  & &  & & &\\  
2.8   &   $T_{2,N}$-Asym & 0.000 & 0.228& 0.658 &&0.652  & 0.980 & 0.996\\
              &   $T_{,N}$-Boot &0.124  &0.668 & 0.904 && 0.944 & 0.990   &0.998 \\ 
              & $T_{N}$-Boot & 0.462& 0.838 & 0.934 && 0.822 &  0.946& 0.986 \\
              &  &  &  & &  & & &\\
3.0    &   $T_{2,N}$-Asym & 0.000 & 0.314& 0.750 && 0.744 &  0.988& 0.998\\
              &   $T_{2,N}$-Boot &0.186 &0.754 &0.908 && 0.966 & 1.000 & 1.000\\
              & $T_{N}$-Boot & 0.530& 0.844 & 0.942 &&  0.834& 0.956 & 0.976 \\                                
         \hline\hline
\end{tabular}
\caption{Empirical power of the tests for the equality of two covariance functions, based on the  statistic $T_{p,N}$-Asym (asymptotic approximation), $T_{p,N}$-Boot and $T_N$-Boot (bootstrap approximation), using two ($p=2$) FPC's, for non-Gaussian data. The curves in the first sample were generated according to (\ref{eq:NG-sim}) while the curves in the second sample were generated according to a scaled version of (\ref{eq:NG-sim}), i.e., $X_2(t)=\gamma X_1(t)$, $t \in {\cal I}$.}
\label{tab:emppower-cov}
\end{table}

\begin{table}[]
\begin{tabular}{|cccc|ccccccc|}
\hline
\hline
                   &              &         &  &  & NG& &&&BB &  \\
       $n_1$ &$ n_2$ & $p$ & Test-Stat. & $\alpha=0.01$ & $\alpha=0.05$ & $\alpha=0.10$ & & $\alpha=0.01$ & $\alpha=0.05$ & $\alpha=0.10$\\
\hline
25 & 25 & 2  &  $S_{p,N}^{(1)}$-Asym & 0.016 &0.062 &0.148 && 0.016 & 0.062 & 0.132\\
&        &    & $S_{p,N}^{(1)}$-Boot & 0.012& 0.052 &0.122 && 0.006 & 0.050 & 0.110 \\ 
     &        &    &   $S_{p,N}^{(2)}$-Asym &0.012 &0.086 & 0.142 && 0.018 & 0.060 & 0.124\\
     &        &    &   $S_{p,N}^{(2)}$-Boot &0.016 &0.042 & 0.084&& 0.018 & 0.060  & 0.132\\
              &        &    &   $S_{N}$-Boot & 0.012& 0.052& 0.122&& 0.016 & 0.064  & 0.136\\
         25 & 25 & 3  &  $S_{p,N}^{(1)}$-Asym &0.020 & 0.078& 0.120&& 0.014& 0.052 & 0.114\\
         &        &    &       $S_{p,N}^{(1)}$-Boot &0.016 &0.042 &0.082 && 0.014 & 0.052 & 0.100\\
     &        &    & $S_{p,N}^{(2)}$-Asym &0.012 &0.048 &0.108 && 0.006 & 0.046 & 0.100\\
     &        &    & $S_{p,N}^{(2)}$-Boot &0.008 &0.036 &0.072 && 0.018 & 0.052 & 0.120\\
              &        &    &   $S_{N}$-Boot &0.016 &0.042 & 0.082&& 0.018 & 0.056  & 0.118\\
              &  &  &  & &  & &&&&\\
              50 & 50 & 2  &  $S_{p,N}^{(1)}$-Asym &0.018 &0.060 &0.114 && 0.008 & 0.066 & 0.130\\
               &        &    &    $S_{p,N}^{(1)}$-Boot &0.008 &0.036 &0.082 && 0.006 & 0.044 & 0.086\\
      &       &    &  $S_{p,N}^{(2)}$-Asym &0.008 & 0.056& 0.114&& 0.008 & 0.068 & 0.114\\
                      &        &    & $S_{p,N}^{(2)}$-Boot &0.010 &0.048 &0.112 && 0.008 & 0.048 & 0.098 \\ 
          &        &    &   $S_{N}$-Boot &0.010 &0.048 &0.112 && 0.008 & 0.046  & 0.100\\
          &  &  &  & &  & &&&&\\
         50 & 50 & 3  & $S_{p,N}^{(1)}$-Asym &0.026 &0.088 &0.134 && 0.016 & 0.070 & 0.116\\
         &        &    &   $S_{p,N}^{(1)}$-Boot &0.010 &0.036 &0.090 && 0.012 & 0.034 & 0.070\\ 
     &        &    &   $S_{p,N}^{(2)}$-Asym &0.010 &0.044 &0.100 && 0.024 & 0.068 & 0.124\\                 
          &        &    &    $S_{p,N}^{(2)}$-Boot & 0.016&0.054 &0.110& &0.060 & 0.052 & 0.102\\
          &        &    &   $S_{N}$-Boot & 0.016&0.054 &0.110 && 0.008 & 0.056  & 0.106\\
                        &  &  &  & &  & &&&&\\
              100 & 100 & 2  &  $S_{p,N}^{(1)}$-Asym &0.014 &0.064 &0.120  && 0.006 & 0.038 & 0.098\\
              &        &    &    $S_{p,N}^{(1)}$-Boot &0.006 &0.042 &0.086 && 0.010 & 0.050 & 0.108\\
      &       &    &  $S_{p,N}^{(2)}$-Asym &0.014 &0.068 &0.122 && 0.002 & 0.042 & 0.098\\               
          &        &    &  $S_{p,N}^{(2)}$-Boot & 0.006& 0.058&0.096 && 0.006 & 0.044 & 0.122\\
          &        &    &   $S_{N}$-Boot &0.006 & 0.058& 0.096&& 0.006 & 0.042  & 0.120\\
          &  &  &  & &  & &&&&\\
         100 & 100 & 3  & $S_{p,N}^{(1)}$-Asym &0.014  &0.062 &0.108 && 0.016 & 0.062 & 0.122\\ 
         &        &    &   $S_{p,N}^{(1)}$-Boot &0.006 &0.046 &0.096 && 0.014 & 0.066 & 0.120\\
     &        &    &   $S_{p,N}^{(2)}$-Asym &0.014 &0.062 &0.110 && 0.016 & 0.076 & 0.136\\
                            &        &    &   $S_{p,N}^{(2)}$-Boot  &0.010 &0.048 &0.092 && 0.016 & 0.068 & 0.112 \\
          &        &    &   $S_{N}$-Boot  &0.010 &0.048 &0.092 && 0.016 & 0.070  & 0.112\\
\hline\hline
\end{tabular}
\caption{Empirical size of the tests for the equality of two mean functions, based on the statistics $S_{p,N}^{(1)}$-Asym and $S_{p,N}^{(2)}$-Asym (asymptotic approximations) and $S_{p,N}^{(1)}$-Boot, $S_{p,N}^{(2)}$-Boot and $S_{N}$-Boot (bootstrap approximations), using two ($p=2$) and three ($p=3$) FPC's, both for Gaussian and non-Gaussian data. The curves in each sample were generated according to Brownian Bridges for Gaussian data (BB)   and according to (\ref{eq:NG-sim}) for non-Gaussian  data (NG).}
\label{tab:empsize-mean}
\end{table}


\begin{table}[]
\centering
\begin{tabular}{|cc|ccccccc|}
\hline
\hline
                          &  &  & N=M=25& &&&N=M=50 &  \\
       $\delta$ & Test-Stat. & $\alpha=0.01$ & $\alpha=0.05$ & $\alpha=0.10$ & & $\alpha=0.01$ & $\alpha=0.05$ & $\alpha=0.10$\\
\hline
1.0 & $S_{2,N}^{(1)}$-Asym & 0.088 & 0.186& 0.304 & &0.142  & 0.314 &0.450 \\
        & $S_{2,N}^{(1)}$-Boot &0.052 &0.160 &0.260 && 0.122 & 0.308 & 0.422 \\ 
    &   $S_{2,N}^{(2)}$-Asym & 0.114 & 0.264& 0.354 &&0.188  & 0.410 &0.522 \\
              &   $S_{2,N}^{(2)}$-Boot & 0.100&0.226 &0.352  &&0.218  & 0.402  & 0.516\\
     &   $S_{N}$-Boot &0.112 &0.270 &0.392  && 0.264  & 0.470  & 0.610\\
              &  &  &  & &  & & &\\              
1.2 & $S_{2,N}^{(1)}$-Asym & 0.138 & 0.280& 0.404 & &0.198  & 0.422 &0.550 \\
        & $S_{2,N}^{(1)}$-Boot &0.076 &0.232 &0.328 && 0.162 & 0.372 &  0.524\\ 
    &   $S_{2,N}^{(2)}$-Asym & 0.168 & 0.332& 0.458 &&0.296  & 0.524 &0.648 \\
              &   $S_{2,N}^{(2)}$-Boot &0.178 &0.338 & 0.440 &&  0.262& 0.524  &0.656 \\
     &   $S_{N}$-Boot &0.202 &0.382 & 0.500 &&0.336  &  0.630 & 0.762 \\
              &  &  &  & &  & & &\\
1.4 & $S_{2,N}^{(1)}$-Asym & 0.170 & 0.356& 0.466 & &0.360  & 0.564 &0.676 \\
        & $S_{2,N}^{(1)}$-Boot &0.118 &0.264 &0.378 && 0.284 & 0.510  & 0.654 \\ 
    &   $S_{2,N}^{(2)}$-Asym & 0.238 & 0.424& 0.530 &&0.430  & 0.676 &0.750 \\
              &   $S_{2,N}^{(2)}$-Boot & 0.214&0.416 & 0.546 && 0.432  &  0.662 & 0.758 \\
     &   $S_{N}$-Boot & 0.246&0.496 &0.616  && 0.536  &  0.776 & 0.878\\ 
                 &  &  &  & &  & & &\\
1.6 & $S_{2,N}^{(1)}$-Asym & 0.262 & 0.448& 0.568 & &0.454  & 0.650 &0.756 \\
        & $S_{2,N}^{(1)}$-Boot &0.134 &0.350 & 0.484&&0.408  &0.638  &0.758  \\ 
    &   $S_{2,N}^{(2)}$-Asym & 0.296 & 0.516& 0.640 &&0.558  & 0.740 &0.830 \\
              &   $S_{2,N}^{(2)}$-Boot &0.284 &0.516 & 0.632 &&0.558  & 0.766  &0.868 \\
     &   $S_{N}$-Boot &0.358 &0.614 &0.728  &&0.686  & 0.888  &0.940 \\ 
           &  &  &  & &  & & &\\
1.8 & $S_{2,N}^{(1)}$-Asym & 0.302 & 0.558& 0.662 & &0.572  & 0.772 &0.852 \\
        & $S_{2,N}^{(1)}$-Boot & 0.196& 0.440& 0.578&& 0.520 & 0.754 & 0.842 \\ 
    &   $S_{2,N}^{(2)}$-Asym & 0.404 & 0.634& 0.722 &&0.672  & 0.858 &0.918 \\
              &   $S_{2,N}^{(2)}$-Boot &0.380 & 0.610&0.708  && 0.686 & 0.866  & 0.914\\
     &   $S_{N}$-Boot &0.458 &0.698 & 0.820 &&0.814  &  0.942 &0.982 \\ 
          &  &  &  & &  & & &\\
2.0 & $S_{2,N}^{(1)}$-Asym & 0.380 & 0.574& 0.690 & &0.668  & 0.824 &0.880 \\
        & $S_{2,N}^{(1)}$-Boot & 0.286&0.512 &0.666 &&0.642  &0.832  &  0.906\\ 
    &   $S_{2,N}^{(2)}$-Asym & 0.434 & 0.658& 0.776 &&0.728  & 0.888 &0.936 \\
              &   $S_{2,N}^{(2)}$-Boot &0.458 & 0.680& 0.794 &&0.792  & 0.914  &0.940 \\
     &   $S_{N}$-Boot &0.576 &0.798 &0.880  &&0.900  & 0.974  &0.994 \\                                                                  
         \hline\hline
\end{tabular}
\caption{Empirical power of the tests for the equality of two mean functions, based on the statistics $S_{p,N}^{(1)}$-Asym and $S_{p,N}^{(2)}$-Asym (asymptotic approximations) and $S_{p,N}^{(1)}$-Boot, $S_{p,N}^{(2)}$-Boot and $S_{N}$-Boot (bootstrap approximations), using two ($p=2$) FPC's, for non-Gaussian data. The curves in the two samples were generated according to the model $X_i(t)=\mu_i(t)+\epsilon_i(t)$ with $\epsilon_i(t)$ distributed independently according to (\ref{eq:NG-sim}), for  $i \in \{1,2\}$, $t \in {\cal I}$. The mean functions were set equal to $\mu_1(t)=0$ and $\mu_2(t)=\delta$.}
\label{tab:emppower-mean}
\end{table} 

\medskip
\newpage

\begin{table}
\begin{tabular}{|c|ccccc|ccccc|}
\hline
\hline
 &              &         &   Absolute && & &&Relative &&  \\
 $T_N$-Boot & 0.179 & &&&& 0.003 && &&\\
 \hline
 &$ p=4$ & $p=5$ & $p=6$ & $p=7$ & $p=8$ &  $p=5$ & $p=6$ & $p=7$ & $p=8$ & $p=9$\\
\hline
$T_{p,N}$-Asym  & 0.253 & 0.211   & 0.385 & 0.545 & 0.520 & 0.004   & 0.021 & 0.064 & 0.130& 0.121\\
$T_{p,N}$-Boot & 0.187 & 0.152   & 0.315 & 0.481 & 0.460 & 0.001  &  0.005 & 0.016 & 0.072& 0.069\\ 
&  &  &  & &  &&&&& \\
$T_{p,N}^{(G)}$-Asym & 0.090 & 0.038  &  0.058 & 0.020 & 0.009 & 0.000  &  0.000 & 0.000 & 0.000& 0.000\\
$T_{p,N}^{(G)}$-Boot & 0.136 & 0.092  &  0.145 & 0.101 & 0.100 & 0.025   & 0.001 & 0.001 & 0.001& 0.002\\ 
 &  &  &  & &  &&&&& \\
 $f_p$  & 0.940 & 0.958   & 0.974 & 0.982 & 0.989 & 0.845   & 0.912 & 0.949 & 0.974& 0.985\\
\hline\hline
\end{tabular}
\caption{$p$-values of the tests for the equality of covariance functions, based on the statistics $T_{p,N}^{(G)}$-Asym and $T_{p,N}$-Asym (asymptotic approximations), $T_{p,N}^{(G)}$-Boot, $T_{p,N}$-Boot and $T_N$-Boot (bootstrap approximations), applied to absolute (left panel) and relative (right panel) egg-laying curves. The term $f_p$ denotes the fraction of the sample variance explained by the first $p$ FPC's, i.e., $f_p = (\sum_{k=1}^p \hat{\lambda}_k)/(\sum_{k=1}^N \hat{\lambda}_k)$.}
\label{tab:medfly-cov}
\end{table}


\begin{table}
\centering
\begin{tabular}{|c|cccccccc|}
\hline
\hline
&  &  &  & &  &&& \\
$S_{N}$-Boot & 0.011& &  & & & & &\\
&  &  &  & &  &&& \\
\hline
 &$ p=2$ & $p=3$ & $p=4$ & $p=5$ & $p=6$ &  $p=7$ & $p=8$ & $p=9$ \\
\hline
$S_{p,N}^{(1)}$-Asym  &0.021  & 0.029   & 0.056 & 0.099 & 0.154 &  0.054  & 0.025 &0.040 \\
$S_{p,N}^{(1)}$-Boot & 0.020& 0.035&0.067 &0.108 & 0.172& 0.070& 0.035 & 0.051 \\
&  &  &  & &  & & & \\
$S_{p,N}^{(2)}$-Asym &0.007 &0.008 &  0.009& 0.009&0.010 &0.010 &0.010  & 0.010 \\ 
$S_{p,N}^{(2)}$-Boot & 0.012&0.011 &0.011  & 0.011& 0.012& 0.011& 0.011&  0.011\\ 
 &  &  &  & &  &&& \\
 $f_p$  &0.837 &0.899 &0.939 & 0.958& 0.973&0.982 &0.989 &0.994  \\
\hline\hline
\end{tabular}
\caption{ $p$-values of the tests for the equality of mean functions, based on the statistics $S_{p,N}^{(1)}$-Asym and $S_{p,N}^{(2)}$-Asym (asymptotic approximations), $S_{p,N}^{(1)}$-Boot, $S_{p,N}^{(2)}$-Boot and $S_N$-Boot (bootstrap approximations), applied to absolute egg-laying curves. The term $f_p$ denotes the fraction of the sample variance explained by the first $p$ FPC's, i.e., $f_p = (\sum_{k=1}^p \hat{\lambda}_k)/(\sum_{k=1}^N \hat{\lambda}_k)$.}
\label{tab:medfly-mean}
\end{table}

\end{document}